\documentclass[leqno,12pt]{article}
\usepackage{latexsym}
\usepackage{amsmath,calligra}
\usepackage{amssymb}
\usepackage{bbm}
\usepackage{color}
\usepackage[latin1]{inputenc}

\usepackage[authoryear]{natbib}
\usepackage{graphicx}
\usepackage{epsfig}
\usepackage{hyperref}

 0

\newcommand{\C}{\mathbb{C}}

\newcommand{\D}{\mathbb T}
\newcommand{\N}{\mathbb{N}}
\newcommand{\Q}{\mathbb{Q}}
\newcommand{\bea}{\begin{eqnarray}}
\newcommand{\eea}{\end{eqnarray}}
\newcommand{\be}{\begin{equation}}
\newcommand{\ee}{\end{equation}}
\newcommand{\beo}{\begin{equation*}}
\newcommand{\eeo}{\end{equation*}}
\newcommand{\beao}{\begin{eqnarray*}}
\newcommand{\eeao}{\end{eqnarray*}}

\newcommand{\V}{\mathcal{V}}
\newcommand{\U}{\mathcal{U}}
\newcommand{\A}{\mathbf{A}}
\newcommand{\X}{\mathbf{X}}
\newcommand{\M}{\mathcal{M}}
\newcommand{\B}{\mathcal{B}}

\newcommand{\tr}{\mathrm{tr}}

\newcommand{\inte}{\mathrm{Int}\,}
\newcommand{\G}{\mathcal{G}}

\newcommand{\br}{\mathcal{B}}

\newcommand{\Int}{\mathrm{Int}\,}
\newcommand{\Clo}{\mathrm{Clo}\,}
\newtheorem{theorem}{Theorem}[section]
\newtheorem{lemma}[theorem]{Lemma}
\newtheorem{prop}[theorem]{Proposition}

\newtheorem{definition}[theorem]{Definition}

\newtheorem{remark}[theorem]{Remark}
\newtheorem{corollary}[theorem]{Corollary}

\makeatletter\@addtoreset{equation}{section}\makeatother

\newcommand{\ba}{\begin{array}}
\newcommand{\ea}{\end{array}}
\newcommand{\beqohne}{\begin{eqnarray*}}
\newcommand{\eeqohne}{\end{eqnarray*}}
\newcommand{\beohne}{\begin{equation*}}
\newcommand{\eeohne}{\end{equation*}}

\def\3{\ss}

\def\en{\mathbb{N}}

\def\ce{\mathbb{C}}

\newcommand{\beq}{\begin{equation}}
\newcommand{\eeq}{\end{equation}}

\def\proof{\noindent{\bf Proof:}\hskip10pt}

\newcommand{\Betam}{\operatorname{Beta}\,}
\newcommand{\GUE}{\operatorname{GUE}\,}

\def \el{\sur{=}{(d)}}
\def \sur#1#2{\mathrel{\mathop{\kern 0pt#1}\limits^{#2}}}
\newcommand{\momi}{\mathcal{M}_n^{[0,1]}}
\DeclareMathAlphabet{\mathcalligra}{T1}{calligra}{m}{n}

\begin{document}

\title{Large Deviations for Random Matricial Moment Problems}

\author{
{\small Fabrice Gamboa}\\
{\small Universit\'e Paul Sabatier}\\
{\small Institut de Math\'ematiques de Toulouse}\\
{\small 118 route de Narbonne}\\
{\small 31062 Toulouse Cedex 9, France}\\
{\small e-mail: gamboa@math.univ-toulouse.fr}\\ 
\and
{\small Jan Nagel} \\
{\small Technische Universit\"at M\"unchen} \\
{\small Zentrum Mathematik} \\
{\small 85747 Garching, Germany} \\
{\small e-mail: jan.nagel@ma.tum.de}\\
\and
{\small Alain Rouault}\\
{\small Universit\'e Versailles-Saint-Quentin}\\
{\small LMV UMR 8100}\\
{\small 45 Avenue des Etats-Unis}\\
{\small 78035-Versailles Cedex France}\\
{\small e-mail: alain.rouault@math.uvsq.fr}
\and
{\small Jens Wagener} \\
{\small Ruhr-Universit\"at Bochum} \\
{\small  Fakult\"at f\"ur Mathematik} \\
{\small 44780 Bochum, Germany} \\
{\small e-mail: jens.wagener@rub.de}\\
}

\maketitle

\begin{abstract}

We consider the moment space  $\mathcal{M}_n^{K}$ corresponding to $p \times p$  complex matrix
measures defined on  $K$ ($K=[0,1]$ or $K=\D$). We endow this set  with the uniform distribution.
We are mainly interested in large deviations principles (LDP) when $n \rightarrow \infty$. First we fix an integer $k$ and study the 
vector of the first $k$ components of a random element of $\mathcal{M}_n^{K}$. 
We obtain a LDP in the set of  $k$-arrays of $p\times p$ matrices. 
Then we lift a random element of $\mathcal{M}_n^{K}$ into a  random measure and prove a LDP at the level of random measures.
We end with a LDP on Carth\'eodory and Schur random functions. These last functions are 
well connected to the above random measure.
In all these problems, we take advantage of the  so-called canonical moments technique by introducing new (matricial) random variables that are independent and have explicit  distributions.

\end{abstract}


\medskip

\section{Introduction}
\subsection{Preliminary: some notations}
All along this article, $p$ will be a positive integer, and $p=1$ will be referred as the scalar case. 
We denote respectively by  
${\mathcal S}_p(\mathbb C)$ the set of all Hermitian $p\times p$ matrices and by ${\mathcal S}_p^+(\mathbb C)$ the one of all Hermitian nonnegative $p\times p$ matrices. If $A, B \in {\mathcal S}_p(\mathbb C)$ we write $A \leq B$ (resp. $A<B$) if, and only if, $B-A$ is nonnegative (resp. positive) definite. 
This is the so-called Loewner partial order on ${\mathcal S}_p(\mathbb C)$  (see for example \cite{hornj85}).
We recall that every $A\in{\mathcal S}_p^+(\mathbb C) $ has a unique nonnegative square root denoted by $A^{1/2}\in{\mathcal S}_p^+(\mathbb C)$. The set of all $p\times p$ unitary matrices is denoted by ${\mathbb U}(p)$.

Let $K$ be either $[0,1]$ or $\mathbb T: = \{ z\in \mathbb C : |z|= 1\}$. 
A matrix-valued probability measure on $K$  is a measure $\mu$ on $K$ 
with values in 
 ${\mathcal S}_p^+(\mathbb C)$ 
  such that
\[\int_K d\mu = I_p\, ,\]
where $I_p$ is the $p\times p$ identity matrix. We denote by $\mathcal{P}(K)$ the set of all matrix-valued probability measures on $K$. In general, if  $(X,{\mathcal A})$ is a measurable 
space, we denote by   $\mathbb M_1(X)$  the set of all probability measures on $X$. We equip it  with the weak convergence topology. This is the 
coarsest topology such that the mappings $\mu \mapsto \int f(x) d\mu(x)$ are continuous. Here, $f \in {\mathcal C}_b(X)$ (the space of bounded continuous functions on $X$) is arbitrary 
(see \cite{moren2008coimbra}  for completeness).

One of the main objects of interest in our work is, for $n\in\mathbb{N}$, the matricial moment space $\mathcal{M}_n^{K}$ defined by 
\begin{equation}
\label{ptitmomo}
\mathcal{M}_n^{K}:=\left\{\left(\int_K x^j d\mu (x)\right)_{ j=1,\ldots,n},\ \mu \in \mathcal{P}(K)\right\}.
\end{equation}
This is a compact set having a nonempty interior - denoted by $\mbox{Int} {\mathcal M}_n^{K}$ -  (see \cite{destu02} for $K=[0,1]$ and \cite{dewag09} for $K = \D$).
 
\subsection{What is done in this paper?}

The aim of our work is to give a picture of the asymptotic behaviour of the set sequence $(\mathcal{M}_n^{K})$. More precisely, we first equip the set $\mathcal{M}_n^{K}$  with the uniform distribution $\mathbb{P}_{K,n}$. Then, for $k \leq n$, we consider $\mathbb{P}_{K,n,k}$ 
the pushforward probability 
 of $\mathbb{P}_{K,n}$ 
under the projection on  $\mathcal{M}_k^{K}$. 
We study, for fixed $k$, the exponential convergence of $(\mathbb{P}_{K,n,k})_n$  when $n$ goes to infinity. 
The asymptotic behavior of  $(\mathbb{P}_{K,n,k})_n$ was widely studied in the scalar case 
beginning with the  seminal paper of \cite{Changetal} where a central limit theorem (CLT) for  $(\mathbb{P}_{[0,1],n,k})$ is proved. Roughly speaking, 
 $(\mathbb{P}_{[0,1],n,k})_n$ converges to the degenerate distribution concentrated on the $k$ first moments of the non symmetric arcsine law and there are Gaussian fluctuations around this limit.
In the same frame, large deviations are studied in \cite{Gamboa}. 
In these papers, the main ingredient for obtaining asymptotic results is a clever  reparametrization of $\mathcal{M}_n^{[0,1]}$. The new parameters, defined recursively, are  
  the so-called canonical moments 
 (see \cite{DeSt97} for a complete overview). Informally, given the the $k-1$ first moments, the $k$-th canonical moment is  the relative position of the $k$-th moment in the range (interval) of possible $k$-th moments. This allows for fixed $n$, to define a bijection between $\mbox{Int} \momi$ and $(0,1)^n$. The key property is that the pushforward of  
 the rather involved probability measure $\mathbb{P}_{[0,1],n,k}$
  under this mapping is a product measure, i.e. the canonical moments are independent.
This is an old result first showed in \cite{Ski} (a simple proof is given in the first chapter of \cite{DeSt97}).
Moreover,  extensions of the asymptotic results on $(\mathbb{P}_{K,n,k})_n$ at the level process are studied in \cite{DeGa}. Also in the scalar case, and using a suitable cousin reparametrization  (also called canonical moments or Verblunsky coefficients) a CLT  and large deviation  are tackled for $(\mathbb{P}_{\D,n,k})_n$ in \cite{LozEJP}. In this last paper, a step toward a multidimensional setting, that is replacing $[0,1]$ by $[0,1]^d\; (d\geq 1)$, is also done. In a more recent work \cite{denag09} 
extend 
 some of the asymptotic results previously described  
 to the matricial moment problem on $[0,1]$ ($p > 1$). As a matter of fact, by using the right extension of canonical moments proposed and first studied in
\cite{destu02}, it is shown there 
 that a CLT holds. 
As before, the key property is the independence, under the uniform distribution on 
$\mathcal{M}_n^{[0,1]}$, of the matricial canonical moment vector.
Here, we revisit these results and obtain new asymptotic result on $\mathcal{M}_n^{K}$. First, we obtain a CLT when $K=\D$. Further, we show large deviations principles (LDP) in both cases, $K=[0,1]$ and $K=\D$. These LDPs are
at level 2, that means that they hold for sequences of distributions of random matricial measures having uniform matricial moments. The main tool 
 is more or less similar as the one used in the scalar case, namely 
  the stochastic independence of the matricial canonical moment. Nevertheless, the matricial case appears to be more technical and due to non commutativity needs more care. Moreover, thanks to the general invariance Proposition \ref{newprop} the complex case ($K=\D$) is tackled by using a polar decomposition 
   argument. 

Besides, it is well known that the truncated trigonometrical problem is connected to two problems of functional analysis on the disc: the so-called Carath\'eodory and Schur problems, respectively. Let us explain the setting in the scalar case, although our results will be in the general matrix case. An analytic function, $F$, on $\mathbb D := \{ z\in \mathbb C : |z| < 1\}$ is called a Carath\'eodory function iff $F(0)= 1$ and $\Re\!\ F(z) > 0$ for all $z \in \mathbb D$. Let ${\mathcal C}_1$ be the set composed by all these functions. An analytic function $f$ on $\mathbb D$ is called a Schur function iff $\sup_{z\in \mathbb D} |f(z)| \leq 1$. Let ${\mathfrak S}_1$ be the set of all Schur functions. The correspondence 
\begin{equation}
\label{1.3.1}
F(z) = \frac{1+zf(z)}{1-zf(z)} \ \ , \ \ f(z) = \frac{1}{z}\!\ \frac{F(z) -1}{F(z) +1}
\end{equation}
is one-one  between ${\mathcal C}_1$ and ${\mathfrak S}_1$.
Any $F \in {\mathcal C}_1$ has a representation
\begin{equation}
\label{1.3.30}F(z) = \int_{\D}\frac{e^{i \theta}+z}{e^{i \theta}-z} d\mu(\theta)\end{equation}
for a unique probability measure $\mu$ on $\mathbb T$ (Herglotz representation theorem).
The Taylor expansion of $F$ is
\begin{equation}
\label{1.3.25}
F(z) = 1 + 2\sum_1^\infty c_n(F) z^n\end{equation}
where  the $c_n$'s are the conjugate  moments of $\mu$, i.e. 
\[c_n(F) =  \int_{\D} e^{-i n\theta} d\mu(\theta) = \bar{\gamma}_n \,.\]
The classical Carath\'eodory problem is to find 
  $F\in \mathcal {C}_1$ such that the first $n$ Taylor coefficients coincide
with given numbers $c_1, \ldots, c_{n}$. It is clearly equivalent to the truncated moment problem.
The Taylor expansion of $f$ is
\begin{equation}
\label{1.3.42a}
f(z) = \sum_0^\infty s_n(f) z^n\,.\end{equation}
 
The Schur problem is to find a Schur function $f(z)$ such that the first $n$ Taylor coefficients coincide
with given numbers $s_0, \ldots, s_{n-1}$. The set \[\mathcalligra{S}_n := \{(s_0(f), \cdots, s_{n-1}(f)) ; f \in {\mathfrak S}_1\}\]
is a compact subset of $\mathbb C^n$.
In the general matrix case, we will study the impact of uniform sampling on the space of Taylor coefficients of these functions. These results are new, even in the scalar case.

One of the main objects
of random matrix theory is to obtain asymptotic results in the limit of large size. Here, on the contrary, the size $p$ of matrices is fixed but
 the dimension  $n$ of the array of matrices  tends to infinity. At first insight, these two topics are very distinct.
Nevertheless, even in the case $p=1$, there is a connection between the random moment problem and the random matrix theory, as described
in \cite{gamrou2010}. Let us formulate it shortly in the generic situation. The  spectral measure of the pair consisting of a  $n \times n$
matrix (unitary or Hermitian) and a fixed vector  is a  discrete measure. It can be  described either by its locations ($n$ points) and its weights, or by a convenient array of its moments. 
When the matrix is random, both representations have  remarkable distributions, and the asymptotical behaviour can be considered from two points of view.
If now we fix $p$ orthonormal vectors instead of only one,
we obtain a random matricial spectral measure and we may consider the array of its (matricial) moments.
This asymptotics will be treated in a forthcoming paper.

The paper is organized as follows.  Section 2 is devoted to the case $K=[0,1]$. It begins with useful definitions and properties around LDPs and ends with the main result on level 2 LDP (Theorem \ref{ldpint}). 
Section \ref{smad} is devoted to the case $K=\D$. We first show a CLT (Theorem \ref{dkonv} and Corollary \ref{dkonv2} ) and then turn to large deviation results (Corollaries \ref{ldpkancirc} and \ref{ldp1}, Theorem \ref{ldp2}). 
In Section 4, we establish a LDP for random Carath\'eodory functions and random Schur functions, respectively (Theorem \ref{thLDPf}).
All technical proofs are postponed to Section \ref{sproo}.
   
\section{Matrix measures on $[0,1]$}
\label{sma01}
Here, we will work on $K=[0,1]$ and the set defined in (\ref{ptitmomo}) is
\begin{align} \label{mom}
\mathcal{M}_n^{[0,1]} := \left \{ \mathbf{S}_n = (S_1,\ldots ,S_n)\ |\ S_j := \int_0^1 x^j d\mu (x),\ j=1,\ldots,n \ \mu \in \mathcal{P}([0,1]) \right \} \subset (\mathcal{S}_p^+(\C))^n,
\end{align}
 The moment space $\mathcal{M}_n^{[0,1]}$ is a compact subset of $(\mathcal{S}_p^+(\C))^n$ with nonempty interior (\cite{destu02}). Therefore the uniform distribution $\mathcal{U}(\mathcal{M}_n^{[0,1]})$ is well defined by the density
\begin{align}
\left(\int_{\mathcal{M}_n^{[0,1]}} dS_1 \dots dS_n\right)^{-1} I\{ \mathbf{S}_n \in \mathcal{M}_n^{[0,1]} \}
\end{align}
with respect to $dS_1\cdots dS_n$ where, if $S= (s_{ij})_{i,j =1}^n$ 
\begin{align} \label{intops}dS= \prod_{i\leq j\leq n} ds_{ij}^\Re\prod_{i<j\leq n}ds_{ij}^\Im\, ,\end{align}
where for $s\in\C$, $s:=s^\Re+i s^\Im$  is the standard decomposition of $s$ in real and imaginary parts. 
 The main tool to study random moments $\mathbf{S}_n \sim \mathcal{U}(\mathcal{M}_n^{[0,1]})$ are the canonical moments which are introduced in the next section.

\subsection{Canonical moments for matrix measures on $[0,1]$}

For a moment vector $(S_1,\ldots ,S_n) \in \momi$ we build
 the block Hankel matrices

\begin{equation} \label{hankel1}
\underline{H}_{2m}:= \left( \begin{array}{ccc}
              S_0 & \cdots & S_m \\
              \vdots &      & \vdots  \\
              S_m & \dots & S_{2m}
           \end{array}
        \right)
~~~~~\overline{H}_{2m}:= \left( \begin{array}{ccc}
       S_1-S_2 & \cdots & S_m-S_{m+1}  \\
              \vdots &      & \vdots  \\
        S_m-S_{m+1} & \dots &S_{2m-1}-S_{2m}
           \end{array}
        \right)
\end{equation}
and
\begin{equation} \label{hankel2}
\underline{H}_{2m+1}:= \left( \begin{array}{ccc}
              S_1 & \cdots & S_{m+1}   \\
              \vdots &      & \vdots  \\
              S_{m+1} & \dots & S_{2m+1}
           \end{array}
        \right)
~~~~~\overline{H}_{2m+1}:= \left( \begin{array}{ccc}
       S_0-S_1 & \cdots & S_m-S_{m+1}  \\
              \vdots &      & \vdots  \\
        S_m-S_{m+1} & \dots &S_{2m}-S_{2m+1}
           \end{array}
        \right).   \\
\end{equation}
\cite{destu02} showed that
the point $(S_1, \ldots , S_n)$ is in $\mbox{Int} \momi$ 
  if, and only if, the matrices $\underline{H}_n$ and $\overline{H}_n$ are both positive definite.

For 
 $(S_1, \ldots , S_n) \in \mbox{Int}(\momi)$  we define
\begin{eqnarray*}
 \underline{h}^{*}_{2m} &:=& (S_{m+1},\cdots,S_{2m})  \\
  \underline{h}^*_{2m-1} &:=& (S_{m},\cdots,S_{2m-1})  \\
 \bar{h}^*_{2m} &:=& (S_{m}-S_{m+1},\cdots,S_{2m-1}-S_{2m})  \\
 \bar{h}^*_{2m-1} &:=& (S_{m}-S_{m+1},\cdots,S_{2m-2}-S_{2m-1})
\end{eqnarray*}
and consider the $p \times p$ matrices
\begin{eqnarray} \label{defs-}
S^-_{n+1}& :=& \underline{h}^*_n \underline{H}^{-1}_{n-1}
\underline{h}_n,~~~n\ge 1~, \\
 \label{defs+}
S^+_{n+1} &:=& S_n - \bar{h}^*_n \bar{H}^{-1}_{n-1} \bar{h}_n, ~~~n\ge 2~,
\end{eqnarray}
(for the sake of completeness we also define $S^-_1 = 0$ and
 $S_1^{+}=I_p$, $S_2^+=S_1)$. Note that $S^-_{n+1}$ and $S^+_{n+1}$
 are continuous functions of $(S_1, \ldots , S_n)$
 and that   $S^-_{n} < S_n < S^+_{n}$  if and only if $(S_1, \ldots , S_{n})
 \in \mbox{Int} \momi$. These preliminary notations allow to introduce the canonical moments of a matrix measure on $[0,1]$.
 
\medskip

\begin{definition}\label{def:kani}
For
 $\mathbf{S}_n=(S_1, \ldots , S_n) \in \inte \momi$ we define the  
 canonical moments by
\begin{align} \label{defkani}
U_k = (S_k^{+} - S_k^{-})^{-1/2} (S_k - S_k^{-}) (S_k^{+} - S_k^{-})^{-1/2},\qquad k=1,\ldots,n
\,.
\end{align}
\end{definition}
\medskip

It is clear that each $U_k \in {\mathcal S}_p(\mathbb C)$ and satisfies
 $0_p < U_k < I_p$. Therefore 
 we can define a mapping 
\begin{align} \label{defphi}
\begin{array}{c} \varphi^{(n)} : \mbox{Int} \momi \longrightarrow (0_p,I_p)^n ,\\
														\varphi^{(n)}(\mathbf{S}_n) = \mathbf{U}_n =  (U_1,\ldots,U_n)\,. \end{array}
\end{align}
By equation \eqref{defkani}, the ordinary moments can be recursively calculated from the canonical moments and the mapping $\varphi^{(n)}$ is one-to-one. Now consider a random vector of moments $\mathbf{S}_n \sim \mathcal{U}(\momi)$, then $\mathbf{S}_n \in \mbox{Int} \momi$ almost surely. 
\cite{denag09} showed that the corresponding canonical moments $\mathbf{U}_n = \varphi^{(n)}(\mathbf{S}_n)$ are independent 
and that $U_k\in {\mathcal S}_p^+(\mathbb C)$ follows a complex  matricial 
 distribution $\Betam\!\!_p(p(n-k+1), p(n-k+1))$ where for $a,b > p-1$ the 
distribution $\Betam\!\!_p(a,b)$ has the density (with respect to $dX$)
\begin{align} \label{betamult}
\B_p(a,b)^{-1} (\det X)^{a -p} (\det (I_p-X))^{b -p}
\end{align}
[see \cite{khatri1965} or \cite{piljou1971}]. The normalizing constant
$\B_p(a,b)$ is defined by 
\begin{align}
\B_p(a,b) := \frac{\Gamma_p(a) \Gamma_p(b)}{\Gamma_p(a+b)}, \qquad a,b> p-1 \ .
\end{align}
Here $\Gamma_p(a)$ denotes the complex multivariate 
 Gamma function
\begin{align*}
\Gamma_p(a) := \pi ^{p(p-1)/2} \prod_{i=1}^p \Gamma (a-i+1), \qquad a > p-1.
\end{align*}
The matricial Beta distribution is one of the three main distributions of complex Hermitian matrices, together with the Gaussian unitary ensemble $\GUE\!\!_p$ having the density
\begin{align} \label{defgue}
( 2  \pi^p)^{-p/2} e^{\displaystyle -\tr \tfrac{1}{2} X^2}
\end{align}
and the complex Wishart distribution $W_p (a)$ 
 with density

\begin{align} \label{defwish}
\Gamma_p(a)^{-1} (\det X)^{a-p} e^{\displaystyle -\tr X},\qquad a> p-1.
\end{align}
We refer to \cite{mehta04} and \cite{forrester2010} for more on these distributions.  
The following result shows that the Wishart distribution and the Gaussian distribution appear as  weak limits of the matricial Beta distribution 
when the parameters tend to infinity.

\medskip

\begin{theorem} \label{betakonv} Let $(a_n)_n$ be a sequence of positive parameters such that $\lim_{n\rightarrow\infty} a_n = \infty$.
\begin{itemize}
\item[(i)]
If $X_n \sim \Betam\!_p(a_n,a_n)$, 
then
\begin{align*}
\sqrt{8a_n}\!\ (X_n - \tfrac{1}{2} I_p) \xrightarrow[n \rightarrow \infty ]{\mathcal{D}} \GUE\!_p\,.
\end{align*} 
\item[(ii)] Let $c>p-1$.
If $X_n \sim \Betam\!_p(c,a_n)$ then 
\begin{align*}
a_n X_n \xrightarrow[n \rightarrow \infty ]{\mathcal{D}} W_p(c)\,.
\end{align*}
\end{itemize}
\end{theorem}

\medskip

The first statement shows that the centered rescaled canonical moments converge in distribution to the $\GUE\!_p$. This is the keystone to obtain a CLT in \cite{denag09}. Notice also, that this implies that the sequence $(X_n)$
converges in probability towards $\tfrac{1}{2} I_p$.
The second statement will play an important role in the study of matrix measures on $\D$.

\medskip

\subsection{Large deviations}

To make this paper self contained let us first recall what is a LDP. For more on LDP we refer to \cite{demboz98}. Let $(u_n)_n$
be an increasing positive sequence of real numbers going
to infinity with $n$.

\begin{definition}\label{dldp}
Let $U$ be a Hausdorff topological space and $\br(U)$ its Borel $\sigma$-field. 
We say that a sequence $(Q_{n})_n$ of probability measures on $(U, \br(U))$ 
  {\rm satisfies a LDP with speed $(u_n)$
and rate function $I : U \rightarrow [0, \infty]$ } if:

\begin{itemize}
\item[i)] $I$ is lower semicontinuous.
\item[ii)] For any measurable set $A$ of $U$:
$$-I(\inte A)\leq
\liminf_{n\rightarrow\infty}u_n^{-1}\log Q_{n}(A)\leq
\limsup_{n\rightarrow\infty}u_n^{-1}\log Q_{n}(A)\leq
-I({\Clo A}),$$
where $I(A)=\inf_{\xi\in A}I(\xi)$ and $\Clo A$ is the closure of $A$.
\end{itemize}
If we omit to give the speed it means that $u_n=n$.
We say that the rate function $I$ is {\rm good} if its level sets
$\{x\in U:\; I(x)\leq a\}$ are compact for any $a\geq 0$.
More generally, a sequence of $U$-valued random variables
is said to satisfy a
LDP if their distributions satisfy a LDP.
\end{definition}
We will need the following well known large deviation result
(see e.g. \cite{demboz98} chapter 4 p. 126 and 130).\\

{\bf Contraction principle}. {\it Assume that $(Q_{n})_n$ satisfies a
LDP on $(U,\br(U))$ with good rate function $I$ and speed $(u_n)$. Let
$T$ be a continuous mapping from $U$ to another Hausdorff topological space $V$.
Then $Q_{n}\circ T^{-1}$ satisfies a LDP on $(V,\br(V))$
with  speed $(u_n)$ and good rate function
\beo
I'(y)=\inf_{x:T(x)=y}I(x),\;\;\;(y\in V).
\eeo}
The so-called cross entropy (or Kullback information) plays an important role in the interpretation of some of our results, for the sake of completeness we recall its definition.
\\

{\bf Kullback Information}.  Let $P$ and $Q$ be probability distributions on $(U,\br(U))$. The Kullback information of $P$ with respect to 
$Q$ is
$${\mathcal{K}}(P ; Q):=
\left\{
\begin{array}{ll}
{\displaystyle \int\log\frac{dP}{dQ} dP,\;\;}&\mbox{ if } P\ll Q \mbox{ and }\log\frac{dP}{dQ}\in L^1(P) \\
\infty &\mbox{ otherwise.}
\end{array}
\right. 
$$
Our first result is a LDP for matricial beta distributions. For the case where the matrix dimension tends to infinity, various LDPs can be found in the literature, see for example \cite{HiaiPetz}. 
Here we are intersted in the case of fixed dimension and growing parameters.

\begin{theorem}\label{ldpbeta}
Let $a_0,a>0$ and $c>p-1$. Further set, for $n\geq 1$, $a_n: = a_0+an$. 

\begin{itemize}
\item[(i)] Let $B_n\sim \Betam\!_p(a_n,a_n)$. Then $B_n$ satisfies a LDP with good rate function
\begin{align}
\mathcal{I}^{(1)}_{B}(B) = 
\begin{cases}
\displaystyle -a \log \det (B-B^2) - 2 ap \log 2,\;\;\;\;&\mbox{if }\;0_p < B < I_p,\\
      \infty\;\;&\mbox{ otherwise.}
\end{cases}
\end{align} 
\item[(ii)] Let $B_n\sim \Betam\!_p(c,a_n)$. Then $B_n$ satisfies a LDP with 
good rate function
\be
\mathcal{I}^{(2)}_B(B)=
\begin{cases}
\displaystyle
-a\log\det(I_p-B),\;\;\;\;&\mbox{if }\;0_p < B < I_p,\\
     \infty\;\;&\mbox{ otherwise.}
     \end{cases}
\ee
\end{itemize}
\end{theorem}

\medskip

\begin{remark}
For the sake of simplicity
we show a LDP only for very special sequences of parameters. This is enough to obtain our further results. However, the result holds for arbitrary sequences $a_n\nearrow\infty$. 
\end{remark}

As a consequence of the last theorem, a LDP for the random matricial vector $\mathbf{U}_k^{(n)} = (U_1,\ldots ,U_k)$ of the first $k$ canonical moments associated to a random matricial vector $\mathbf{S}_n$  uniformly drawn holds.  
Indeed, as mentioned before, the components of $\mathbf{U}_k^{(n)} = (U_1,\ldots ,U_k)$ are independent, so that we obtain:

\begin{corollary} \label{ldpkan}
Let $\mathbf{S}_n \sim \mathcal{U}(\momi)$ and for $k$ fixed, let $\mathbf{U}_k^{(n)}$ 
denote the projection of $\mathbf{U}_n = \varphi^{(n)}(\mathbf{S}_n)$ onto the first $k$ coordinates. 
Then  the sequence $\left(\mathbf{U}_k^{(n)}\right)_n$ satisfies a LDP in $({\mathcal S}_p^+(\ce))^k$ with 
good rate function
\begin{align} 
\mathcal{I}_{\mathbf{U}}(\mathbf{U}_k) = 
\begin{cases}\displaystyle- \sum_{i=1}^k p \log \det (U_i-U_i^2) - 2k p^2  \log 2,\;\;\;\;&\mbox{if}\;\mathbf{U}_k \in (0_p,I_p)^k,\\
     \infty\;\;&\mbox{ otherwise.}
     \end{cases}
\end{align} 
\end{corollary}

\medskip

Obviously the rate function $\mathcal{I}_{\mathbf{U}}$ achieves its minimum value $0$ at $\mathbf{U}_k = (\frac{1}{2}I_p,\ldots ,\frac{1}{2}I_p)$ that appears as discussed before for general sequences of matricial beta distributed random matrices, see Theorem \ref{betakonv}) as the limit of $\mathbf{U}_k^{(n)}$. 
Notice also that the constant infinite sequence $U_k = \frac{1}{2}I_p$ , $k \geq 1$ is the moment sequence of the matrix arcsine law $\nu_p$ defined by 
\begin{align} \label{defarcsin}
d\nu_1 (x) = \frac{dx}{\pi \sqrt{x(1-x)}} \ \ , \ \ d\nu_p(x) = d\nu_1(x) I_p \ , \ (p>1)\,,
\end{align}
see \cite{denag09}.

Now, the vector of  ordinary moments $(S_1,\ldots ,S_k)$ is a continuous function of the canonical moment vector $\mathbf{U}_k^{(n)}$. So we obtain the following Corollary from Corollary \ref{ldpkan} by a simple application of the contraction principle and the identity 
\begin{align} \label{detid}
\det (S_{k+1}^+ - S_{k+1}^-) = \det \prod_{i=1}^k U_i (I_p - U_i)
\end{align}
(see \cite{destu02}).

\begin{corollary} \label{ldpgew}
Let $\mathbf{S}_n \sim \mathcal{U}(\momi)$ and for $k<n$ let $\mathbf{S}_k^{(n)}$ denote the projection of $\mathbf{S}_n$ onto the first $k$ coordinates. Then  $\mathbf{S}_k^{(n)}$ satisfies a LDP with 
good rate function
\begin{align} 
\mathcal{I}_{\mathbf{S}}(\mathbf{S}_k) = 
\begin{cases}- p \log \det (S_{k+1}^+ - S_{k+1}^-) - 2k p^2  \log 2,
\;\;\;\;&\mbox{if}\;\mathbf{S}_k \in \inte \momi,\\
     \infty\;\;&\mbox{ otherwise.}
     \end{cases}
\end{align} 
\end{corollary}
\medskip

 We end this section with a LDP for random matrix measures  on $[0,1]$. For this purpose,  for every $n$ let $\mathbb{P}_n$ denote 
 any probability measure on  $\mathcal{P}([0,1])$ 
such that the pushforward by the mapping
\[\mu \in \mathcal{P}([0,1]) \mapsto \mathbf{S}_n(\mu) = (S_1 (\mu),\ldots ,S_n (\mu)) \in \momi\]
is $\mathcal{U}(\momi)$.


\medskip

\begin{theorem} \label{ldpint}
The sequence 
 $(\mathbb{P}_n)_n$ satisfies a LDP in $\mathbb M_1({\mathcal P}([0,1]))$ with 
 good rate function
\begin{equation}
\mathcal{I}_{[0,1]}(\mu) = 
\begin{cases}
\displaystyle-p \int_0^1 \log \det W(x)\, d\nu_1(x), \;\;\;\;&\mbox{if}\;\nu_1\{ \det W = 0\} = 0,\\
      \infty\;\;&\mbox{ otherwise.}
\end{cases}
\end{equation}
where 
$d\mu(x) = W(x) d\nu_p(x) + d\mu^s(x)$  
is the Lebesgue decomposition\footnote{see \cite{robertson}  on Lebesgue decomposition for matricial measures} of $\mu$ with respect to 
 $\nu_p$ as matricial measures on $[0,1]$ ($\nu_1$ and $\nu_p$ are the arcsine measures defined by \eqref{defarcsin}).
 \end{theorem}
\medskip

\begin{remark}
\label{rem1}
\begin{enumerate}
\item When $p=1$ (scalar case) the rate function is also
\begin{equation}
\label{K1}
\mathcal{I}_{[0,1]}(\mu) ={\mathcal K}(\nu_1 ; \mu)\,. 
\end{equation} 
The matricial case has also an interpretation in terms of cross-entropy which we hope to address in a future work.
\item
A cousin 
 result 
 of Theorem \ref{ldpint} holds 
 in the frame of real matrix measures. In this case the constant $p$ in the rate function is replaced by $\frac{p+1}{2}$. 
 All arguments remain essentially unchanged and we refer to \cite{denag09} for the underlying 
 results on real matrix valued random moments and the corresponding canonical moments. 
 
\item
 From Theorem \ref{ldpint} and Corollary \ref{ldpgew} together with the contraction principle 
 one easily obtains the following identity of rate functions.
 For $\mathbf{S}_k = (S_1,\dots, S_k) \in\inte \momi$ we have
\begin{align}
\mathcal{I}_{\mathbf{S}}(\mathbf{S}_k) = - p \log \det (S_{k+1}^+ - S_{k+1}^-) - 2k p^2  \log 2
= \inf_{\mathcal{D}(\mathbf{S}_k)} -p \int_0^1 \log \det W(x) d\nu_1(x),
\end{align}
where
\begin{align}
\mathcal{D}(\mathbf{S}_k) = \left \{ \mu \in \mathcal{P}([0,1])\ |\ \int_0^1 x^j d\mu (x) = S_j, \ j = 1,\dots ,k \right \}
\end{align}
and $W$ is defined as in Theorem \ref{ldpint}.
\end{enumerate}
\end{remark}

\bigskip

\section{Matrix measures on $\D$: the trigonometric case}

\label{smad}
In this section, we consider the space $\mathcal{P}(\D)$ of  matrix-valued probability measures on the unit circle $\D$.
In what follows $\Gamma_j$ denotes the $j-$th trigonometric moment of a matrix measure $\mu\in\mathcal{P}(\D)$, that is
\be
\Gamma_j=\Gamma_j(\mu)=\int_{-\pi}^{\pi}e^{ij\theta}d\mu(\theta)
\ee
and for $n\in\N$ and $p\geq1$ the set defined in (\ref{ptitmomo}) is 
\begin{align}
\mathcal{M}_n^{\D}&:=\left\{(\Gamma_1,\dots,\Gamma_n)|\ \Gamma_j=\Gamma_j(\mu), \mu\in\mathcal{P}(\D)\right\}\subset\left(\C^{p\times p}\right)^n.
\end{align}
Unlike to moments of matrix measures on $[0,1]$, the moment $\Gamma_j$ is no more Hermitian.
Therefore we use the following Lebesgue measure
on $\C^{p\times p}$. 
For $X\in\C^{p\times p}$ define
\be\label{integraloperator}
dX 
=\prod_{1\leq i,j\leq p}d
x_{ij}^\Re d
x_{ij}^\Im
\,.
\ee

\subsection{Canonical moments on $\D$}

As in the above section we use a  notion of canonical moments to study  $\mathcal{M}_n^{\D}$. 
First,  for $(\Gamma_1,\dots,\Gamma_n)\in\mathcal{M}_n^{\D}$,  we build the block Toeplitz matrix
\be\label{toepmat}
T_n:=\left(\Gamma_{i-j}\right)_{i,j=0,\dots,n}.
\ee
\cite{dewag09} showed that  $(\Gamma_1,\dots,\Gamma_n)\in \Int \mathcal{M}_n^{\D}$ 
if and only if $T_n>0$. Therefore this interior is non empty. 
 Furthermore they proved that 
for $(\Gamma_1,\dots,\Gamma_n)\in\inte \mathcal{M}_n^{\D}$ the range of the moment $\Gamma_{n+1}$ is  the set
\be
K_n=\left\{W\in\C^{p\times p}\ |\ L_n^{-1/2}(W-M_n)R_n^{-1/2}=U, \ UU^{\ast}\leq I_p\right\},
\ee
where the matrices $L_n$, $R_n$ and $M_n$ are defined  by
\begin{align}\label{lmrm}
L_n&:=\left[I_p-\left(\Gamma_1, \dots ,\Gamma_n\right)T_{n-1}^{-1}\left(\Gamma_1, \dots ,\Gamma_n\right)^{\ast}\right],\\
R_n&:=\left[I_p-\left(\Gamma_{-n}, \dots ,\Gamma_{-1}\right)T_{n-1}^{-1}\left(\Gamma_{-n}, \dots ,\Gamma_{-1}\right)^{\ast}\right]\label{rm},\\
\label{mm}
M_n&:=\left(\Gamma_1, \dots ,\Gamma_n\right)T_{n-1}^{-1}\left(\Gamma_{-n}, \dots ,\Gamma_{-1}\right)^{\ast},
\end{align}
respectively. In this frame, canonical moments are defined by normalizing the moments 
 in the following way.

\begin{definition}\label{kancirc}
For $(\Gamma_1,\dots,\Gamma_n)\in \inte \mathcal{M}_n^{\D}$ we define the canonical moments $A_j$, $j=1, \cdots, n$ setting 
\be
A_1:=\Gamma_1,\quad A_{j}:=L_{j-1}^{-1/2}(\Gamma_{j}-M_{j-1})R_{j-1}^{-1/2}\ \ (j=2, \cdots, n)\,.
\ee
\end{definition}

The canonical moments of a matrix measure always lie in the set
\be
\mathbbm{D}_p=\{U\in\C^{p\times p}\ |\ UU^{\ast}\leq I_p\}
\ee
and coincide with the well known Verblunsky coefficients appearing in the Szeg\"{o} recursion of orthonormal matrix polynomials (see e.g. \cite{simon05} Section 2.13). They are connected to the trigonometric moments by a one-to-one mapping $\psi^{(n)}:\inte\mathcal{M}_n^{\D}\rightarrow\inte\mathbbm{D}_p^n$ recursively defined by Definition \ref{kancirc}.

\medskip

We now state a Taylor expansion of the inverse of the mapping $\psi^{(n)}$. Here and in the following $\|M\|$ always denotes the Frobenius norm of the complex entries matrix $M$, that is 
$$\|M\|:=\tr(M^*M)^{1/2}.$$

\begin{lemma}\label{taylorentwicklung}
Let $n\in\N^{+}$ and $\A_n=(A_1,\dots,A_n)\in \inte\mathbbm{D}_p^n$. The mapping $(\psi^{(n)})^{-1}:\A_n\mapsto\X_n=(\Gamma_1,\dots,\Gamma_n)$ induced by the definition of canonical moments has an order one Taylor expansion at  $0$. Namely,
\be
\X_n=\A_n+o(\|\A_n\|).
\ee
\end{lemma}

In the following this Taylor expansion will be used to derive results concerning trigonometric moments from results obtained for canonical moments.

\subsection{Weak convergence in the trigonometrical case}

As in the real case we define a uniform distribution $\U(\M_n^{\D})$ on $\mathcal{M}_n^{\D}$ by the density
\be
\left(\int_{\M_n^{\D}}d\Gamma_1\dots d\Gamma_n\right)^{-1}I\left\{\X_n\in\M_n^{\D}\right\},
\ee
now with respect to the measure (\ref{integraloperator}).We first state a result on the distribution of the canonical moments when the corresponding trigonometric moments are uniformly distributed.

\begin{lemma}\label{vertkanmom}
Let $\X_n\sim\U(\M_n^{\D})$ and $\A_n=(A_1,\dots,A_n)=\psi^{(n)}(\X_n)\in (\mathbb D_p)^n$ denote the corresponding vector of canonical moments. Then 
$A_1, \dots , A_n$ are independent 
 and for $k=1, \dots, n$, $A_k$ has density
\be\label{dichtea}
\frac{1}{c_k^{(n)}}\det\left(I_p-A_k^*A_k\right)^{2p(n-k)}
\ee
with respect to (\ref{integraloperator}), where $c_k^{(n)}$ is a normalizing constant.
\end{lemma}

We now establish a relation between the Hermitian random matrices from Section \ref{sma01} and matricial random variables 
without symmetry condition:

\begin{theorem}
\label{newtheo}
If $A_k$ is a random matrix with density (\ref{dichtea}), then
\begin{equation}
\label{eloi}
A_k \el VB_k^{1/2}
\end{equation}
where $V$ and $B_k$ are independent, $V$ is Haar distributed in $\mathbb U(p)$ and $B_k$ follows a multivariate complex Beta distribution $\Betam\!_{\!p}(p,2p(n-k)+p)$ (see \ref{betamult}).
\end{theorem}

The previous theorem is a particular case of the following general variable change result. It is quite natural  and useful in other asymptotical  problems involving random complex matrices. Similar 
arguments have been used recently by \cite{fischmann1} to generate matrices of the Ginibre ensemble.

\begin{prop}
\label{newprop}
Let $M$ be  a $p\times p$ random matrix with complex entries whose density with respect to (\ref{integraloperator}) is 
$f(x_1^2(M), \cdots, x_p^2(M))$ where $x_1(M), \cdots, x_p(M)$ are the (positive) singular values, and $f$ is a symmetric function. Then, the random matrices 
$H= M^*M$ and $U=\left(M^*M\right)^{-1/2}M$ are independent, $U$ is Haar distributed in $\mathbb U(p)$ and 
the density of $H\in {\mathcal S}_p^+(\mathbb C)$ with respect to (\ref{intops}) is proportional to $f(\lambda_1(H), \cdots, \lambda_p(H))$ 
where $\lambda_1(H), \cdots, \lambda_p(H)$ are the eigenvalues of $H$.  
\end{prop}

We are now in the position to give our first limit theorem in the trigonometrical case.

\begin{theorem}\label{dkonv}
Let $\X_n\sim\U(\M_n^{\D})$, $\A_n=\psi^{(n)}(\X_n)$ and $\A_n^k$ denote the projection onto the first $k$ coordinates ($k$ is fixed). 
Then for $n\rightarrow\infty$ the weak convergence
\be
\sqrt{2pn}\A_n^k\overset{\mathcal{D}}{\longrightarrow}\G_k
\ee
holds, where $\G_k=(G_1,\dots,G_k)$ and $G_1,\dots,G_k$ are complex iid random matrices of the Ginibre complex ensemble
(see \cite{gigi}), or, in other words,	having density
\begin{equation}
\label{gauss}
g(G)=\pi^{-p^2}\exp{(-\|G\|^2)}
\end{equation}
with respect to (\ref{integraloperator}).
\end{theorem}

As a consequence, using the Taylor expansion of Lemma \ref{taylorentwicklung} and the $\delta$-method (see for example \cite{vava}),
 we obtain a weak convergence theorem for the  rescaled random trigonometric moments. This is the subject of the next corollary.

\begin{corollary}\label{dkonv2}
Let $\X_n\sim\U(\M_n^{\D})$ and $\X_n^k$ denote the projection onto the first $k$ coordinates ($k$ is fixed). Then when $n\rightarrow\infty$ 
\be
\sqrt{2pn}\X_n^k\overset{\mathcal{D}}{\longrightarrow}\G_k,
\ee
(here $\G_k$ is  as in Theorem \ref{dkonv}). 
\end{corollary}


\subsection{Large deviations in the trigonometrical case}

Our final results concern LDPs for random moments and matrix measures on the unit circle. The large deviations in the scalar trigonometrical case are due to  \cite{LozEJP} Theorems 4.2 and 4.4. Nevertheless, in that paper, there was a mistake in the computation of the Jacobian. A power $2$ is missing.

The proof of the next Corollary  follows directly from part (ii) of Theorem \ref{ldpbeta} (applying  the contraction principle). We again use the equality  $A_k \el VB_k^{1/2}$, where $B_k\sim \Betam\!_p(p,2p(n-k)+p)$ and $V$ is Haar distributed on the unitary group. By Lemma \ref{vertkanmom} the canonical moments are independent, giving the final form of the rate function.

\begin{corollary}\label{ldpkancirc}
Let $\X_n\sim\U(\M_n^{\D})$, $\A_n=\psi^{(n)}(\X_n)$ and $\A_n^k$ denote the projection onto the first $k$ coordinates ($k$ is fixed). Then $\A_n^k$ satisfies a LDP with 
good rate function
\be
\mathcal{I}_{\A}(\mathbf{Z})=
\mathcal{I}_{\A}(Z_1,\dots,Z_k)=\begin{cases}\displaystyle -2p\sum_{i=1}^k\log\det\left(I_p-Z_i^*Z_i\right),\;\;\;\;&\mbox{if}\;\mathbf{Z}\in \inte\mathbbm{D}_p^k,\\
      \infty\;\;&\mbox{ otherwise.}
\end{cases}
\ee
\end{corollary}

Another application of the contraction principle for the mapping $\psi^{(n)}$ yields the following LDP for the trigonometric moments.

\begin{corollary}\label{ldp1}
Let $\X_n\sim\U(\M_n^{\D})$ and $\X_n^k$ denote the projection onto the first $k$ coordinates ($k$ is fixed). Then $\X_n^k$ satisfies a LDP 
with 
good rate function
\be
\mathcal{I}_{\Gamma}(\X)=\mathcal{I}_{\Gamma}(\Gamma_1,\dots,\Gamma_k)=
\begin{cases}\displaystyle -2p\log{\frac{\det(T_k)}{\det(T_{k-1})}},\;\;\;\;&\mbox{if}\;\X\in \inte\M_k^{\D},\\
      \infty\;\;&\mbox{ otherwise.}
\end{cases}
\ee
 Here,  $T_k$ denotes the block Toeplitz matrix (\ref{toepmat}) defined by $(\Gamma_1,\dots,\Gamma_k)$.
\end{corollary}

Finally we state a LDP for a sequence of random matrix measures on $\D$. 
For every $n$, let $\Q_n$ denote 
a probability measure on the set $\mathcal{P}(\D)$ such that 
the pushforward by the mapping
\[\mu \in \mathcal{P}(\D) \mapsto \mathbf{X}_n(\mu) = (\Gamma_1 (\mu),\ldots ,\Gamma_n (\mu)) \in \M_n^{\D}\]
is $\mathcal{U}(\M_n^{\D})$.

\begin{theorem}\label{ldp2}
The sequence 
 $(\Q_n)_n$ satisfies a LDP in $\mathbb M_1\left({\mathcal P}(\D)\right)$ 
 with 
good rate function
\be
\mathcal{I}_{\D}(\mu)=
\begin{cases}\displaystyle-\frac{p}{\pi}\int_{\D}\log{\det(W(\theta))}d\theta,\;\;\;\;&\mbox{if}\;\det W(\theta)\neq0\;\mbox{a.e.},\\
     \infty\;\;&\mbox{ otherwise,}
\end{cases}
\ee
 where
 $d\mu(\theta)=W(\theta)\frac{d\theta}{2\pi}+d\mu^s(\theta)$ 
is the Lebesgue decomposition of $\mu$ with respect to $\frac{d\theta}{2\pi} I_p$ as matricial measures 
on $\D$.
\end{theorem}

The proof is very similar to that one of Theorem \ref{ldpint} and therefore omitted. 

\begin{remark}
\begin{enumerate}
\item For $p=1$ the rate function is also
\begin{equation}
\label{K2} 
\mathcal{I}_{\D}(\mu) = 2 {\mathcal K}\left(\frac{d\theta}{2\pi}  ; \mu\right)\,.
\end{equation}
It is the content of Theorem 4.4 in \cite{LozEJP} but a factor $2$ was missing  in that paper, owing to  a mistake in the Jacobian (7.2).
\item
As in Remark \ref{rem1} we see, 
from Theorem \ref{ldp2} and Corollary \ref{ldp1} together with the contraction principle, 
 the following identity of rate functions.
For $\X_k=(\Gamma_1,\dots,\Gamma_k)\in\inte\mathcal{M}_k^{\D}$ we have
\be
\mathcal{I}_{\Gamma}(\X_k)=-2p\log\frac{\det(T_k)}{\det(T_{k-1})}=\inf_{\mathcal{C}(\X_k)}-\frac{p}{\pi}\int_{\D}\log\det(W(\theta))d\theta,
\ee
where
\be
\mathcal{C}(\X_k)=\left\{\mu\in\mathcal{P}(\D)\ \left|\ \int_{-\pi}^{\pi}e^{ij\theta}d\mu(\theta)=\Gamma_j,\ j=1,\dots,k\right.\right\}
\ee
and $W$ is defined as in Theorem \ref{ldp2}.
\end{enumerate}
\end{remark}

\section{Application: Random Carath\'{e}odory and Schur matrix functions} 

In the above Theorem \ref{ldp2}, we studied a family of random measures. 
Since the truncated trigonometrical moment problem is closely connected to the Carath\'eodory problem, which is itself connected to the Schur problem, it may be natural to look at the corresponding random functions. 
 In this section we study the impact of uniform  sampling on the space of Taylor coefficients of these functions. We first give the framework, which can be seen in \cite{damanik2008analytic} or \cite{dubovoj1992matricial}  and then we give our results. It seems to be new, even in the scalar case.  

\subsection{Carath\'eodory and Schur matrix-valued functions}
As before,
let $p$ be a given positive integer. By a $\mathbb C^{p\times p}$-valued Carath\'eodory matrix function  $F(z)$, one means a $p\times p$ matrix-valued function which is holomorphic in $\mathbb D$, has a nonnegative real part there
\[F^\Re (z) \equiv \frac{1}{2} (F(z) + F(z)^*) \geq 0, \ \ \ z \in \mathbb D\,,\]
and such that $F(0) = I_p$.
We use the notation ${\mathcal C}_p$ to designate the class of such $\mathbb C^{p\times p}$-valued Carath\'eodory matrix functions.
We also define the class $\mathfrak{S}_p$  of $\mathbb C^{p\times p}$-matrix valued functions $f$ analytic in $\mathbb D$ and contractive there, i.e. such that $f(z) \in \bar{\mathbb D}_p$ for $z \in \mathbb D$ 
, which are called matrix valued Schur functions. 

 The correspondence
\be
F(z)=(I_p+zf(z))(I_p-zf(z))^{-1}\quad\text{and}\quad f(z)=z^{-1}(F(z)-I_p)(F(z)+I_p)^{-1}
\ee
is one-to-one between ${\mathcal C}_p$ and $\mathfrak S_p$. 
Any $F \in {\mathcal C}_p$ has a representation
\[F(z) = \int_{\D}\frac{e^{i \theta} +z}{e^{i \theta} -z}\ d\mu(e^{i \theta})\,, \ z \in \mathbb D\,,\]
for a unique   $\mu \in {\mathcal P}(\D)$. Any $F \in {\mathcal C}_p$ has a finite radial limit 
$\lim_{r\uparrow 1} F(re^{i\theta}) =: F(e^{i\theta})$ for almost every $\theta$. The corresponding value of $f$ in such a point $e^{i\theta}$ will be denoted by $f(e^{i\theta})$. 
If 
\beo
d\mu(\theta)=W(\theta)\frac{d\theta}{2\pi}+d\mu_s(\theta)
\eeo
is the Lebesgue decomposition of $\mu$ one has the identity
\begin{equation}
\label{1.3.32}
W(\theta)= F^\Re(e^{i\theta})=(I_p-e^{-i\theta}f(e^{i\theta})^*)^{-1}(I_p-f(e^{i\theta})^*f(e^{i\theta}))(I_p-e^{i\theta}f(e^{i\theta}))^{-1}
\end{equation}
a.e. and for a.e. $\theta$, $\det W(\theta) \not= 0$ iff $f(e^{i\theta})^*f(e^{i\theta}) < 1$ (Prop. 3.16 in \cite{damanik2008analytic}).

The Taylor expansion of $F$ is given by
\beo
F(z)= 
I_p + 2\sum_{k=1}^{\infty}C_k(F)z^k,
\eeo
where the coefficients 
 are the conjugate trigonometric moments of the matrix measure $\mu$ associated to $F$, i.e.
\[ C_k(F) = \int_\D e^{-ik\theta} d\mu(\theta) = \Gamma_k^*\,.\]
The classical Carath\'eodory problem is to find 
  $F\in \mathcal {C}_p$ such that the first $n$ Taylor coefficients coincide
with given $p\times p$ matrices $C_1, \ldots, C_n$. It is clearly equivalent to the truncated moment problem.

Each Schur function in $\mathfrak{S}_p$ is associated to a matrix measure $\mu\in {\mathcal P}(\D)$, hence to the sequence of its canonical moments $(A_k)_{k \geq 1}$. For every $j \geq 1$, let $f_j$ be the Schur function corresponding to the shifted sequence
$(A_k)_{k \geq j+1}$, and set $f_0 =f$. From Theorem 3.19 of \cite{damanik2008analytic}, we have the recursive relations:
\begin{eqnarray}
\label{Schurmat1}
f_{k}(z) &=& z^{-1}(B_k^R)^{-1}[f_{k-1}(z) - A_{k}^*][I_p - A_kf_{k-1}(z)]^{-1} B_k^L\,,\\
\label{Schurmat2}
f_k(z) &=& (B^R_{k+1})^{-1}[z f_{k+1}(z)+ A_{k+1}^*][I_p + zA_{k+1}f_{k+1}]^{-1} B_{k+1}^L\,,
\end{eqnarray}
where 
\begin{equation}
B_k^R := [I_p - A_k^*A_k]^{1/2} \ ,\  B_k^L := [I_p - A_kA_k^*]^{1/2}\,.
\end{equation}
The Taylor expansion of $f$ is
\begin{equation}
\label{1.3.42}
f(z) = \sum_0^\infty G_k(f) z^k\,.\end{equation}

The Schur problem is to find a Schur function $f\in \mathfrak S_p$ such that the first $n$ Taylor coefficients coincide
with given numbers $G_0, \ldots, G_{n-1}$.  A solution
exists if and only if the block matrix
\begin{align*}
\begin{pmatrix} G_0& 0& 0& \ldots &0\\ G_1& G_0&  0& \ldots &0\\ G_2& G_1& G_0& \ldots &0\\
&&\ldots&\\ G_{n-1}&G_{n-2}&G_{n-3}&\ldots &G_0
\end{pmatrix}
\end{align*}
is contractive, i.e. if it satisfies $GG^* \leq I_{np}$ (see \cite{dubovoj1992matricial}, Theorem 3.1.1). The set
\[\mathcalligra{S}_n := \{(G_0(f), \cdots, G_{n-1}(f)) ; f \in {\mathfrak S}_p\}\]
is a relatively compact subset of $(\mathbb C^{p\times p})^n$.

In both problems, the system of canonical moments (alias Verblunsky coefficients, alias Schur coefficients) plays a prominent role. 
In  Section \ref{vertkanmom} we saw that the dependence between the moments (hence the $C_k$'s) and the canonical moments is triangular. 
The relation between the Taylor coeffcients of a Schur function and its Schur coeffcients (i.e. the canonical moments of the associated measure)
 is also triangular. We postpone the presentation of this point in
 the proof of Theorem \ref{thLDPf}.   

\subsection{Randomization. Large deviations}

\medskip

 For every $n$ let $\mathbb{P}^c_n$ denote a probability measure on the set ${\mathcal C}_p$ 
  such that 
the pushforward by the mapping
\[F \in \mathcal{C}_p \mapsto \mathbf{C}_n(F) = (C_1(F) ,\ldots ,C_n(F) ) \in \mathcal{M}_n^{\D}\]
is $\mathcal{U}(\mathcal{M}_n^{\D})$.
Let also $\mathbb{P}^s_n$ denote a probability measure on the set ${\mathfrak S}_p$ such that 
the pushforward by the mapping
\[f \in \mathfrak{S}_p \mapsto {\mathbf G}_n(f) := (G_0(f),\ldots ,G_{n-1}(f)) \in \mathcalligra{S}_n\]
is $\mathcal{U}(\mathcalligra{S}_n)$.

One gets the following LDP for matrix valued Carath\'eodory and Schur functions.

\begin{theorem}
\label{thLDPf}
The sequence $(\mathbb{P}^c_n)_n$ satisfies a LDP in $\mathbb M_1({\mathcal C}_p)$ with good rate function
\begin{equation}
\label{LDPMatCara}
\mathcal{I}_p^C(F) =\begin{cases}\displaystyle - \frac{p}{\pi}\int_{\D} \log \det  F^\Re(e^{i\theta})d\theta,\;\;\;\;&\mbox{if}\;\det F^\Re(e^{i\theta})\neq0\;\mbox{a.e.},\\ 
      \infty\;\;&\mbox{ otherwise.}
\end{cases}
\end{equation}
The sequence $(\mathbb{P}^s_n)_n$ satisfies a LDP in $\mathbb M_1(\mathfrak{S}_p)$ with good rate function
\begin{equation}
\label{LDPMatSchur}
\mathcal{I}_p^S(f) =\begin{cases}\displaystyle -\frac{p}{\pi}\int_{\D}\log\det (I_p-f(e^{i\theta})^*f(e^{i\theta}))d\theta,
\;\;\;\;&\mbox{if}\;\det (I_p-f(e^{i\theta})^*f(e^{i\theta}))\neq0\;\mbox{a.e.},\\
      \infty\;\;&\mbox{ otherwise.}
\end{cases}
\end{equation}
\end{theorem}

\begin{remark}
 Behind Theorem \ref{ldp2} and Theorem \ref{thLDPf} (and as will be seen in the proofs), there is a triple identity, which holds true in the generic case: 
\begin{eqnarray}
\nonumber
\sum_n \log \det (I_p - A_n A_n^*) 
&=& \int_{\D} \log \det W(\theta)\frac{d\theta}{2\pi} 
= \int_{\D} \log \det F^\Re(e^{i\theta})\frac{d\theta}{2\pi}\\ 
 &=& \int_{\D}\log\det (I_p-f(e^{i\theta})^*f(e^{i\theta}))\frac{d\theta}{2\pi},
\end{eqnarray}
say 
\[(1) = (2)=(3)=(4).\]
Equality $(1)= (2)$ is  Szeg\"o's Theorem for matrix-valued measures (see Theorem 2.13.5 in \cite{simon05}), and $(1)= (4)$ is the matricial version of Boyd's theorem (see 2.7.7 of \cite{simon05} in the scalar case).
\end{remark}

\section{Proofs}
\label{sproo}

\subsection{Proof of Theorem \ref{betakonv}}

If $X$ is $\Betam\!_p(\alpha, \beta)$ distributed, then
\[X \el \left(W_1 + W_2\right)^{-1/2} W_1 \left(W_1 + W_2\right)^{-1/2}\]
where $W_1 \sim W_p(\alpha)$ and $W_2 \sim W_p(\beta)$ are independent 
and Wishart distributed.

For $(i)$, we choose $\alpha = \beta = a_n$ and observe that
\[ X_n - \tfrac{1}{2} I_p \el \frac{1}{2}\!\ \left(W_1 + W_2\right)^{-1/2} \left[\left(W_1-a_n I_p\right) + \left(a_n I_p - W_2\right)\right]
\left(W_1 + W_2\right)^{-1/2}
\]
then we apply Proposition \ref{LLNCLTW} (i) and (ii).

For $(ii)$, it is enough to take $\alpha = c$ and $\beta= a_n$ and apply Proposition \ref{LLNCLTW} (i). 
\hfill $\Box$

\subsection{Proof of Theorem \ref{ldpbeta}}
We give a proof only for $a_n=an$.\\
To prove (i) let $B_n \sim \Betam\!_p(an, an )$, then again the following equality in distribution holds
\begin{align} \label{betawish}
B_n \el \left( \sum_{i=1}^{2n} W_i \right)^{-1/2} \left( \sum_{i=1}^n W_i\right) \left( \sum_{i=1}^{2n} W_i \right)^{-1/2},
\end{align}
where the random variables 
 are independent and $W_p(a)$ distributed. 
  (see e.g. \cite{piljou1971}). 
  By Proposition \ref{wish1} each component $V_n^{(1)}, V_n^{(2)}$ of the vector
\begin{align*}
\begin{pmatrix} V_n^{(1)} \\ V_n^{(2)} \end{pmatrix}
=
\begin{pmatrix}
							\frac{1}{n} \sum_{i=1}^n W_i \\ \frac{1}{n} \sum_{i=n+1}^{2n} W_i \end{pmatrix}
\end{align*}
satisfies a LDP with good rate function $\Lambda^\star$ given by (\ref{rf1}).

The independence of the random variables $W_i$ now yields a LDP for  $(V_n^{(1)}, V_n^{(2)})$ with good rate function
$\Lambda^\star (X) + \Lambda^\star (Y)$. By the contraction principle and equality \eqref{betawish} the random variable $B_n$ satisfies a LDP on $(0_p, I_p)$ with good rate function
\begin{align*}
\mathcal{I}(Z) =& \inf_{\mathcal{Z}} \left( \Lambda^\star (X) + \Lambda^\star (Y) \right) \\
								=& \inf_{\mathcal{Z}} 
												\left( \tr (X+Y) - a \log \det (XY) - 2pa + 2pa \log a \right),
\end{align*}
where the infimum is taken over the set 
\begin{align*}
\mathcal{Z} = \left \{ (X,Y) \in \mathcal{S}_p^+ (\ce)^2\ |\ Z = (X + Y)^{-1/2}X(X + Y)^{-1/2} \right \}\,.
\end{align*}
On $\mathcal{Z}$ we have $\det(XY) = \det (Z(I_p-Z)\det(X+Y)^2)$ and we can write the rate function as
\[\mathcal{I}(Z) = -a\log\det \left(Z(I_p-Z)\right) - 2pa + 2pa \log a + \inf_{\mathcal{Z}} 
											\left( \tr (X+Y) - 2a \log \det (X+Y)\right)\,.\]
Appealing to (\ref{hineq}) with $L = (2a)^{-1}(X+Y)$, we see that 
\begin{align*}
\mathcal{I}(Z) = -a \log \det (Z(I_p-Z)) - 2pa \log 2.
\end{align*}

\medskip
 
To prove (ii) let $B_n \sim \Betam\!_p(c, an )$. Then we have
\[B_n \el \left(\frac{X}{n} + \frac{1}{n}\sum_{i=1}^nW_i\right)^{-1/2} \frac{X}{n} \left(\frac{X}{n} + \frac{1}{n}\sum_{i=1}^nW_i\right)^{-1/2}\]
where $X\sim W_p(c)$,  $(W_i)_{i=1,\ldots,n}$ are iid $W_p(a)$ distributed and $X$ and $(W_i)_{i=1,\ldots,n}$ are independent. 
By Propositions \ref{wish1} and   \ref{newlem}, we get for $\left(\frac{X}{n}, \frac{1}{n}\sum_{i=1}^nW_i\right)$ a LDP with rate function the sum of rate functions and by the contraction principle, we get a LDP with rate function 
\[\mathcal{I}(Z)=\inf_{\mathcal Z} (\tr X + \tr Y -a\log\det Y -ap+ap\log a),\] 
where ${\mathcal Z}$ is as in the proof of Theorem \ref{ldpbeta} (i). 
On ${\mathcal Z}$ we have $\det(Y) = \det(X+Y) \det (I_p - Z)$, hence 
\[\tr X + \tr Y -a\log\det Y = \tr(X+Y) - a\log \det (X+Y) -a\log\det(I_p - Z)\]
and the infimum is achieved for $(X+Y) = aI_p$ by (\ref{hineq}). This completes the proof.\hfill $\Box$

\subsection{Proof of Theorem \ref{ldpint}}

We follow
 here the proof given in \cite{Gamboa} concerning the scalar case. 
Let  $\widetilde{\mathbb{P}}_n$ be the probability measure on the infinite dimensional moment space 
\begin{align*}
\mathcal{M}_{\infty}^{[0,1]} = \left \{ \mathbf{S} = (S_1,S_2,\ldots )\ |\ S_j = \int_0^1 x^j d\mu (x),\ \mu \in \mathcal{P}([0,1]) \right \}
\end{align*}
 induced by the bijection $\mathbf{S} \mapsto \mu_{\mathbf{S}}$. Now if $\Pi_k^{\infty}$ denotes the canonical projection 
$\mathcal{M}_{\infty}^{[0,1]} \rightarrow \mathcal{M}_{k}^{[0,1]}$, then the measure 
$\tilde{\mathbb{P}}_n \circ \left(\Pi_k^{\infty}\right)^{-1}$ is the law of $\mathbf{S}_k^{(n)}$. Therefore, Corollary \ref{ldpgew} yields a LDP for the sequence $\left( \tilde{\mathbb{P}}_n \circ \left(\Pi_k^{\infty}\right)^{-1} \right)_n$ with speed $n$ and good rate function
\begin{align*}
\widetilde{\mathcal{I}}_k(\mathbf{S}_k) = - p \log \det (S_{k+1}^+ - S_{k+1}^-) - 2k p^2  \log 2.
\end{align*}
By Dawson-G\"artner's Theorem (see \cite{demboz98}) the sequence $\widetilde{\mathbb{P}}_n$ satisfies a LDP with good rate function
\begin{align*}
\widetilde{\mathcal{I}}(\mathbf{S}) = \sup_{k \in \en} \widetilde{\mathcal{I}}_k(\mathbf{S}_k).
\end{align*}
It remains to calculate the right hand side of the last equality, which is given by 
\begin{align*}
\sup_{k \in \en} - p \log \left( 4^{pk} \det (S_{k+1}^+ - S_{k+1}^-) \right) .
\end{align*}

Let $\mu$ denote a matrix measure corresponding to the sequence $\mathbf{S}_k$ and let $\tilde{\mu}$ denote the image measure on $[-1,1]$ obtained from 
$\mu$ by the affine transformation $x \mapsto 2(x-\tfrac{1}{2})$. 
Since canonical moments are invariant under affine transformations, i.e., $U_i(\mu)=U_i(\tilde{\mu})$ (see for example \cite{denag09}, Lemma 3.1), we have
\begin{align*}
\det (S_{k+1}^+(\mu) - S_{k+1}^-(\mu)) = \prod_{i=1}^k \det(U_i(\mu) - U_i^2(\mu))  = \prod_{i=1}^k \det(U_i(\tilde{\mu}) - U_i^2(\tilde{\mu})),
\end{align*}
where the first identity is again \eqref{detid}. Now denote by $\mu_C$ the symmetric matrix measure on  $\D$ associated with $\tilde{\mu}$, that is
\begin{align}
\int_{-1}^1 f(x) d\tilde{\mu}(x) = \int_{-\pi}^{\pi} f(\cos(\theta))d\mu_C(\theta).
\end{align}
The canonical moments $U_i(\tilde{\mu})$ are related to the canonical moments $A_i(\mu_C)$ by the relation (see \cite{dewag09})
\begin{align*}
U_i(\tilde{\mu}) = \frac{1}{2} (A_i(\mu_C) + I_p).
\end{align*}
This gives for the range
\begin{align*}
\det (S_{k+1}^+(\mu) - S_{k+1}^-(\mu)) = \prod_{i=1}^k 4^{-p} \det(I_p - A_i(\mu_C)^2).
\end{align*}
Since $0\leq\det(I_p - A_i(\mu_C)^2)\leq1$, the sequence $\widetilde{\mathcal{I}}_k(\mathbf{S}_k)$ is increasing in $k$ which yields
$$
\sup_{k \in \en} - p \log \left( 4^{pk} \det (S_{k+1}^+ - S_{k+1}^-) \right)=\lim_{k \rightarrow \infty} - p \log \left( \prod_{i=1}^k \det(I_p - A_i(\mu_C)^2) \right).
$$
Then the Szeg\"o's Theorem for Matrix-Valued Measures (Theorem 2.13.5 in \cite{simon05}) yields
\begin{align*}
\widetilde{\mathcal{I}}(\mathbf{S}(\mu)) =& \lim_{k \rightarrow \infty} - p \log \left( \prod_{i=1}^k \det(I_p - A_n(\mu_C)^2) \right)\\
								=& -\frac{p}{2\pi} \int_{-\pi}^{\pi} \log \det W(\theta) d\theta ,
\end{align*}
where $d\mu_C(\theta) = W(\theta) \frac{d\theta}{2\pi} + d\mu_S$ is the Lebesgue decomposition of $\mu_C$. 
 Since $\mu_C$ is symmetric, $W$ is an even function
\begin{align*}
\widetilde{\mathcal{I}}(\mathbf{S(\mu)}) =& -\frac{p}{\pi} \int_{0}^{\pi} \log \det W(\theta) d\theta
\end{align*}
which, after projection on $[0,1]$ yields 
\begin{align*}
\widetilde{\mathcal{I}}(\mathbf{S(\mu)}) =& -\frac{p}{\pi} \int_{0}^{1} \log \det V(x) \frac{dx}{\sqrt{x(1-x)}}
\end{align*}
where $V(x) = W(\arccos(2x-1))$ is the Radon-Nikodym derivative of $\mu$ with respect to the arcsine matricial measure. The result follows from the contraction principle and the continuity of the mapping $\mathbf{S} \mapsto \mu_{\mathbf{S}}$. \hfill $\Box$

\medskip

\medskip

\subsection{Proof of Lemma \ref{taylorentwicklung}}

First we recall the  notion of Fr\'{e}chet differentiability (see for example \cite{oldcar}).\\
Let $\U$ be an open subset of a complex Banach space $X$ and $\Phi$ a continuous map from $\mathcal U$ to a complex Banach space $Y$. The map $\Phi$ is called differentiable  at $U\in \U$, if there exists a bounded linear operator $L$ from $X$ to $Y$ such that
\[\lim_{V\rightarrow 0} \frac{\Vert \Phi(U+V) - \Phi(U) - LV\Vert}{\Vert V\Vert} =0\,.\]
 We denote $L$ by $D\Phi(U)$ 
  and call it differential of $\Phi$ at $U$.\\
For this notion of differentiability we have the following rules :
\begin{itemize}
\item $(chain-rule)$ Let $Z$ be a Banach space, $\mathcal V$ be an open subset of $Y$ and  $\Psi:\V\rightarrow Z$ be a continuous mapping from $\mathcal V$ to $Z$. If $\Phi(U)\in\V$, if $\Phi$ is differentiable at $U$ and if $\Psi$ is differentiable at $\Phi(U)$ then $\Psi\circ\Phi$ is differentiable at $U$ and
\be
D(\Psi\circ\Phi)(U) = D\Psi(\Phi(U))\circ D\Phi(U)\,.
\ee
\item $(product-rule)$ If we have a multiplicative structure on $Y$ 
 and
if $\Phi$ and $\Psi$ are continuous maps from $\U$ to $Y$, both differentiable  at $U_0$ 
 then the map  $\Phi\Psi : U \mapsto \Phi(U)\cdot \Psi(U)$ is differentiable at $U_0$ and for every $V$
 \be
 D(\Phi\Psi)(U_0)V = \left[D\Phi(U_0)V\right]\cdot\Psi(U_0) + \Phi(U_0) \cdot\left[D\Psi(U_0)V\right] .
 \ee
\end{itemize}
We note that the mapping $M\mapsto M^{1/2}$ is differentiable at $I_p$. Further,  the action of the differential at that point is the multiplication by $\frac{1}{2}$. Theorem \ref{taylorentwicklung} now follows using the above mentioned rules and the following lemma.

\begin{lemma}\label{rekursionl}
Let $(\Gamma_1,\dots,\Gamma_n)\in\Int\mathcal{M}_n^{\D}$. For the matrices $L_n$ and $R_n$ defined in (\ref{lmrm}) and (\ref{rm}), respectively, the following recursions hold
\be
L_n=L_{n-1}^{1/2}\left(I_p-A_nA_n^{\ast}\right)L_{n-1}^{1/2}\quad\text{and}\quad R_n=R_{n-1}^{1/2}\left(I_p-A_n^{\ast}A_n\right)R_{n-1}^{1/2}.
\ee
\end{lemma}

\medskip

{\bf Proof:}
We only show the result for $L_n$. For $R_n$, the proof is left for the reader.\\
Here we use the notation of \cite{dewag09}. Let $\phi_n^L$ and $\phi_n^R$ be the orthonormal matrix polynomials. 
 Using the Szeg\"{o} recursion (compare e.g. \cite{simon05} section 2.13) and the fact that $L_n^{-1/2}$ is Hermitian we obtain
\begin{align*}
I_p&=\langle z\phi_n^L,z\phi_n^L\rangle_L\\
&=\langle L_n^{-1/2}L_{n+1}^{1/2}\phi_{n+1}^L+A_{n+1}\tilde{\phi}_n^R,\ L_n^{-1/2}L_{n+1}^{1/2}\phi_{n+1}^L+A_{n+1}\tilde{\phi}_n^R\rangle_L\\
&=L_n^{-1/2}L_{n+1}L_n^{-1/2}+A_{n+1}A_{n+1}^{\ast}\,.
\end{align*}
Indeed the definition of the inner products directly yields
\beo
\langle\tilde{\phi}^R_n,\tilde{\phi}^R_n\rangle_L=\langle\phi^R_n,\phi^R_n\rangle_R=I_p.
\eeo
The assertion of the Lemma follows.

\hfill $\Box$

In the following we will differentiate mappings from $\C^{np\times p}$ to 
$\C^{p \times p}$. 
We have from the definition of canonical moments
\be\label{darstgamma}
\Gamma_k=L_{k-1}^{1/2}A_kR_{k-1}^{1/2}+M_{k-1} \ \ \ (1 \leq k \leq n)\,,
\ee
where the matrices $L_{k-1}$, $R_{k-1}$ and $M_{k-1}$ are defined in (\ref{lmrm}) to (\ref{mm}). The differentiability of $\A_n\mapsto L_{k-1}^{1/2}A_kR_{k-1}^{1/2}$ at $0_p^{(n)}=(0_p,\dots,0_p)\in\C^{np\times p}$ follows obviously using the product rule. Indeed, first the linear map $\A_n\mapsto A_k$ is obviously differentiable in $0_p^{(n)}$. The action  of the differential is the multiplication by  the map itself. 
 The  differentiability of $\A_n\mapsto L_k$ and $\A_n\mapsto R_k$ can be established using induction on $k$ and Lemma \ref{rekursionl} together with chain and product rules. Again by induction one obtains $L_k(0_p^{(n)})=R_k(0_p^{(n)})= I_p$.  Now the product rule yields, for every $V \in \mathbb C^p$
\begin{align*} 
D(L_{k-1}^{1/2}A_kR_{k-1}^{1/2}) 
(0_p^{(n)})V &= \big[DL_{k-1}^{1/2}(0_p^{(n)})V\big]\cdot A_k(0_p^{(n)})\cdot  
 R_{k-1}^{1/2} 
(0_p^{(n)})+ 
 L_{k-1}^{1/2} (0_p^{(n)})\cdot A_k V 
  R_{k-1}^{1/2}(0_p^{(n)})\\ &+ (L_{k-1}^{1/2} (0_p^{(n)})\cdot A_k(0_p^{(n)})\cdot \big[D R_{k-1}^{1/2}(0_p^{(n)})V\big]  \\
&= A_k V\,.
\end{align*}
It remains to show that $M_{k-1}=o(\|\A_n\|)$ for $k=1,\dots,n$.
It is done by induction with respect to $k$ together with an appeal to the continuity of the inversion at $I_{(k-1)p}$. This yields the conclusion of Lemma \ref{taylorentwicklung}.
\hfill $\Box$

\medskip

\subsection{Proof of Lemma \ref{vertkanmom}}
We have by definition of the canonical moments that $A_k$ depends only on $\Gamma_1,\dots,\Gamma_k$ so that the Jacobian of $\psi^{(n)}$ is the product of the Jacobians of $(\Gamma_1,\dots,\Gamma_k)\mapsto A_k$ ($k=1,\dots,n$). As
\beo
A_k=L_{k-1}^{-1/2}(\Gamma_k-M_{k-1})R_{k-1}^{-1/2}
\eeo
and because $L_{k-1}$, $R_{k-1}$ and $M_{k-1}$ are independent of $\Gamma_k$, Theorem 3.2 from \cite{mathai97} gives the following Jacobian $J_k$ for the mapping $\Gamma_k\mapsto A_k$:
\begin{align*}
J_k&=\det\left(L_{k-1}^{-1/2}\left(L_{k-1}^{-1/2}\right)^{\ast}\right)^{p}\det\left(R_{k-1}^{-1/2}\left(R_{k-1}^{-1/2}\right)^{\ast}\right)^{p}\\
&=\det(L_{k-1})^{-p}\det(R_{k-1})^{-p},
\end{align*}
where the last equality follows because $L_{k-1}$ and $R_{k-1}$ are Hermitian.
From Lemma \ref{rekursionl} we obtain
\beo
\det(L_{k-1})^{-p}\det(R_{k-1})^{-p}=\prod_{j=1}^{k-1}\det(I_p-A_j^*A_j)^{-p}\det(I_p-A_j^{\ast}A_j)^{-p}=\prod_{j=1}^{k-1}\det(I_p-A_jA_j^{\ast})^{-2p}.
\eeo
Consequently, the Jacobian of $\psi^{(n)}$ is the product
\begin{align*}
\prod_{k=1}^n\prod_{j=1}^{k-1}\det\left(I_p-A_j^*A_j\right)^{2p} = \prod_{k=1}^{n-1}\det\left(I_p-A_k^*A_k\right)^{2p(n-k)}
\end{align*}
This yields exactly 
 the assertion of the lemma.
\hfill $\Box$

\subsection{Proof of Proposition \ref{newprop}}
The proof of this proposition uses the following lemma.

\begin{lemma}
\label{newlemma}
Let $A$ be a $p\times p$ matrix of full rank and $A = U H^{1/2}$ its polar decomposition with $H=A^*A\in {\mathcal S}_p(\mathbb C)$ 
and  $U = A (A^*A)^{-1/2}\in \mathbb U(p)$. If $A$ is random and if
\begin{equation}
\label{hyp}
\forall V \in \mathbb U(p) \ \ A \el VA
\end{equation}
then $U$ and $H$ are independent, and $U$ is Haar distributed.
\end{lemma}

\noindent{\bf Proof of Lemma \ref{newlemma}} 

We have for all bounded measurable functions $f_1,f_2$
\begin{eqnarray}
\nonumber
\mathbb E \left(f_1(U) f_2(H)\right) &=&  \mathbb E f_1\left( A (A^*A)^{-1/2}\right) f_2\left((A^*A)\right)\\
\label{1}
&=& \mathbb E f_1\left( VA (A^*A)^{-1/2}\right) f_2\left((A^*A)\right)\\
\label{2}
&=& \int_{\mathbb U(p)} \left[\mathbb E f_1\left( VA (A^*A)^{-1/2}\right) f_2\left((A^*A)\right)\right] d_{Haar}(V)\\
\label{3}
&=& \mathbb E\left(  \left[\int_{\mathbb U(p)} f_1\left( VA (A^*A)^{-1/2}\right)  d_{Haar}(V)\right] f_2\left((A^*A)\right)\right)\\
\label{4}
&=& \mathbb E\left(  \left[\int_{\mathbb U(p)} f_1( V)  d_{Haar}(V)\right] f_2\left((A^*A)\right)\right)\\
\label{conclusion}
&=& \left[\int_{\mathbb U(p)} f_1( V)  d_{Haar}(V)\right] \mathbb E\left(f_2\left((A^*A)\right)\right),
\end{eqnarray}
where in (\ref{1}) we take into account the invariance by left multiplication, in (\ref{2}) the fact that  $V$ is arbitrary in $\mathbb U(p)$, in (\ref{3}) Fubini's theorem, and in (\ref{4}) the invariance of Haar by right multiplication.
\hfill $\Box$

\noindent{\bf Proof of Proposition \ref{newprop}}

The assumption (\ref{hyp}) is trivially verified since $VA$ and $A$ have the same singular values. It remains to determine the distribution of  
 $H= M^*M$. By a simple application of Proposition 4.1.3 of \cite{agz}, we see that the singular values of $M$ 
 have on $(0,\infty)^p$ a joint density proportional to
 \[|\Delta(x_1^2, \cdots, x_p^2)|^2 f (x_1^2, \cdots, x_p^2) (x_1\dots x_p)\]
 where $\Delta$ is the Vandermonde function. 
This implies directly that the eigenvalues of $H$ have on $(0,\infty)^p$ a joint density proportional to
\[|\Delta(\lambda_1, \cdots, \lambda_p)|^2 f (\lambda_1, \cdots, \lambda_p)\,.\]
 Now it is easy to lift to the matrix $H$ by Proposition 4.1.1 of  \cite{agz}.
\hfill $\Box$
\medskip
 
\noindent{\bf Proof of Theorem \ref{newtheo}} 
If $A_k$ has density $f(A_k)$ it fulfills the assumptions of Proposition \ref{newprop}, 
with
\[f(\lambda_1, \cdots, \lambda_p) = \frac{1}{c_k^{(n)}}\prod_{j=1}^p (1- \lambda_j)^{2p(n-k)}\]
and the density of $B_k$ is proportional to
\[\det (I_p - B_k)^{2p(n-k)}\,.\]
This expression fits with (\ref{betamult}) with
$a=p$ and $b= 2p(n-k) + p$. 
\hfill $\Box$

\subsection{Proof of Theorem \ref{dkonv}}

One proof of Theorem \ref{dkonv} directly follows from two applications of Theorem \ref{newtheo} together with Lemma \ref{vertkanmom}, Theorem \ref{betakonv} and the continuous mapping theorem. We give a second proof here.

\subsubsection{Alternative proof: Gaussian approximation}

We  use two clever results. The first one will give a representation of the law of $A_k$.
\begin{theorem}[\cite{collins2005product} Theorem 5.1 or \cite{forrester2009derivation}]
 The top $p\times p$ sub-block of a Haar distributed matrix from $\mathbb U(p+q)$, where $q \geq p$, has a density in $\mathbb D_p$ proportional to
\[A \mapsto \det \left(I_p - A A^*\right)^{q-p}\,.\]
\end{theorem}
The second one is the following "Borel theorem".
\begin{theorem}[\cite{jiang2005maxima}, Corollary 1] 
There exists two $N\times N$ random matrices $\Pi_N = (\pi_{i,j})_{1\leq i, j\leq N}$ and  
$Y_N = ( y_{i,j})_{1\leq i, j\leq N}$ 
 defined on the same probability space such that
\begin{itemize}
\item[i)] $\Pi_N$ is Haar distributed in $\mathbb U(N)$

\item[ii)] all the $y_{i,j}, 1\leq i, j\leq N$ are independent and standard complex gaussian distributed. 

\item[iii)] For $m_N = [N/ (\log N)^2]$
\[\max_{i\leq N, j \leq m_N} |\sqrt{N} \pi_{i,j} - y_{i,j}| \rightarrow 0\]
in probability as $N \rightarrow \infty$.
\end{itemize}
\end{theorem}

From the above notation and Lemma \ref{vertkanmom},  $A_k$ is distributed as the top $p\times p$ sub-block of $\Pi_N$
with $N= 2p(n-k+1)$.  
Up to a change of probability space we have then for $i,j \leq p$
\[\sqrt{2p(n-k+1)}(A_k)_{i,j} - y_{i,j} \rightarrow 0\]
in probability as $n \rightarrow \infty$, which leads easily to the conclusion since $k$ is fixed. 
\hfill $\Box$

\subsection{Proof of Corollary \ref{ldp1}}

By the contraction principle and Corollary \ref{ldpkancirc},  $(\X_n^k)_n$ satisfies a LDP with 
good rate function
\beo
\widetilde{\mathcal{I}}_{\Gamma}(\Gamma_1,\dots,\Gamma_k)=\begin{cases}-2p\sum_{i=1}^k\log{\det(I_p-A_i^*A_i)},&\text{ if }(\Gamma_1,\dots,\Gamma_k)\in \inte\M_k^{\D},\\
\infty&\text{ otherwise},\end{cases}
\eeo
where $(A_1,\dots,A_k)=\psi^{(k)}(\Gamma_1,\dots,\Gamma_k)$. An application of the formula for determinants of block matrices (see for example \cite{hornj85}) yields
\beo
\det(T_k)=\det(T_{k-1})\det(R_k)=\det(T_{k-1})\det(L_k),
\eeo
because $L_k$ and $R_k$ are Schur complements in $T_k$. From Lemma \ref{rekursionl} we obtain
\beo
\det(R_k)=\prod_{i=1}^k\det(I_p-A_i^*A_i)
\eeo
and so
\beo
\sum_{i=1}^k\log{\det(I_p-A_i^*A_i)}=\log{\frac{\det(T_k)}{\det(T_{k-1})}},
\eeo
which is the assertion of Corollary \ref{ldp1}. \hfill $\Box$

\subsection{Proof of Theorem \ref{thLDPf}}

For $(\mathbb{P}^c_n)$ (Carath\'eodory problem), the assertion is a consequence of  Theorem \ref{ldp2}, the contraction principle and (\ref{1.3.32}). Recall the main  point: under $\mathcal{U}(\mathcal{M}_n^{\D})$, the variables $A_1, \cdots, A_n$ are independent, and $A_k$ has a density proportional to 
$\det\left(I_p-A_j^*A_j\right)^{2(n-j)p}$.

For $(\mathbb{P}^s_n)$ (Schur problem), we first remark from (\ref{Schurmat1}) that the mapping $(G_0(f), \cdots, G_{n-1}(f)) \mapsto (A_1, \cdots, A_n)$ is triangular, i.e. that $G_k(f)$ depends only on $A_1, \cdots, A_{k+1}$. Let us give details.  In the scalar case, it is  1.3.48 in \cite{simon05} and we follow the same scheme, up to change due to non commutativity. 
Relation (\ref{Schurmat2}) for $k=0$ implies
\[f(z)(B_1^L)^{-1}[I_p +zA_1f_1(z)] = (B_1^R)^{-1} [zf_1(z) + A_1^*]\,.\]
Identifying the powers of $z^n$ on both sides yields
\begin{eqnarray*}
G_0(f) &=& (B_1^R)^{-1}A_1^* B_1^L\\
G_n(f) &=& (B_1^R)^{-1}G_{n-1}(f_1)B^L_1 - G_0(f) (B_1^L)^{-1}A_1 G_{n-1}(f_1) - \sum_{j=1}^{n-1} G_j(f) (B_1^L)^{-1}A_1 G_{n-1-j}(f_1)
\end{eqnarray*}
Lemma 1.3 in \cite{damanik2008analytic} (see also formula (2.13.52) in \cite{simon05}) says that
\[A_j^*B^L_j = B^R_jA_j^*\]
for every $j \geq 1$ so that we get
$G_0(f) = A_1$ and identifying the powers of $z^n$ on both sides yields:
\begin{eqnarray}
\nonumber
G_0(f) &=& A_1^*\\
G_n(f) &=& (B_1^R)^{-1} 
G_{n-1}(f_1) B_1^L - \sum_{j=0}^{n-1} G_j(f) (B_1^L)^{-1}A_1 G_{n-1-j}(f_1)\ (n\geq 1)\,.
\end{eqnarray}
Induction on $n$ leads to
\begin{eqnarray}
\nonumber
G_n(f) &=& V_n A_{n+1}^* W_n  
\\
&+& \mbox{polynomial in }\ (A_1, A_1^*, \cdots, A_n, A_n^*)\,. 
\end{eqnarray}
where 
\begin{eqnarray*}
V_n = B_1^R
B_2^R
 \cdots B_{n}^R
 \ , \ 
W_n = B^L_n B_{n-1}^L \cdots B_1^L\,.
\end{eqnarray*}
From this relation, we see that, if we froze $A_1, \cdots, A_n$ the Jacobian of the mapping $G_n(f) \mapsto A_{n+1}$ is (Theorem 3.2 of \cite{mathai97})
\[|\det(V_nV_n^*)|^p |\det(W_nW_n^*)|^p = \prod_{k=1}^n [\det (I_p - A_k^*A_k)]^{2p}\,.\]
Like in the proof of Lemma \ref{vertkanmom},  it turns out that the Jacobian of the mapping 
\[(G_0(f), \cdots, G_{n-1}(f)) \mapsto (A_1, \cdots, A_{n})\]
is then
\[\prod_{k=1}^{n-1} \det(I_p - A_k^*A_k)^{2(n-k)}\,.\]
We conclude that
the distribution of 
$(A_1, \cdots, A_n)$ under  $\mathbb{P}^s_n$ is the same as the distribution of $(A_1, \cdots, A_n)$ under $\mathbb{P}^c_n$.  Applying again the contraction principle, we see that 
$(\mathbb{P}^s_n)$ satisfies a LDP with good rate function
\[I_p^s (f) = - \frac{p}{\pi} \int_\D \log \det W(\theta) d\theta\]
where $W$ is related to $\mu$ the underlying matrix measure. To have a rate function depending explicitly on $f$, we go back to the correspondence (\ref{1.3.32}) between $W$ and $f$ so that
\[\log \det W(\theta) = \log \det (I_p-f(e^{i\theta})^*f(e^{i\theta})) -2 \log |\det (I_p-e^{i\theta}f(e^{i\theta}))| \]
 and apply Jensen's formula to the function  
$\det (I_p-zf(z))$. This yields (\ref{LDPMatSchur}).
 \hfill $\Box$

\section{Appendix: some properties of the Wishart distribution}
\label{sAppen}
 
For $a > 0$, the Laplace transform of the complex Wishart distribution $W_p(a)$ is given for $K \in {\mathcal S}_p$ by
 \begin{align}
\Lambda(K) = \log \mbox{E} \left[ e^{\tr (KW)} \right] = - a \log \det (I_p-K)
\end{align}
if $K < I_p$ and infinite otherwise.
From the divisibility of the family of Wishart distributions (indexed by $a$), we deduce the following easy results (law of large numbers and CLT).
\begin{prop}
\label{LLNCLTW}
As $a_n \rightarrow \infty$ we have for $W_n \sim W_p(a_n)$
\begin{enumerate}
\item[(i)]
$\displaystyle\lim_{n\rightarrow\infty} \frac{1}{a_n} W_n =  I_p\ \ \ \hbox{(in probability)}\,,$
\item[(ii)]
$\displaystyle(a_n)^{-1/2} \left(W_n -a_n I_p\right)\overset{\mathcal{D}}{\longrightarrow} \GUE\!_p$
\end{enumerate}
\end{prop}

Since the following large deviations result is not so obvious, we give a proof.

 \begin{prop}
 \label{wish1}
 For fixed $p$ and $a > 0$, if the variables $X_k, k \geq 1$ are independent and  $W_p(a)$ distributed, then $\frac{1}{n}(X_1+ \cdots+X_n)$ satisfies a LDP in 
 ${\mathcal S}_p^+(\C)$ with good rate function 
\begin{align} \label{rf1}
\Lambda^\star (X) =\begin{cases}\displaystyle \tr X - a\log \det X - ap (1 - \log a)\;\;\;\;&\mbox{if}\;\det X > 0,\\
      \infty\;\;&\mbox{ otherwise.}
\end{cases}
\end{align}
 \end{prop}

\proof
The multidimensional Cram\'er theorem gives a LDP with good rate function
\be\label{lambdastar}
\Lambda^\star(X) = \sup_{K\in{\mathcal S}_p(\C)} \tr (KX) - \Lambda(K).
\ee
We first give a non variational  expression of $\Lambda^\star(X)$.

If $\det X =0$, for every $n$ we choose 
 $K_n \in {\mathcal S}_p(\C)$ such that $K_nx= 0$ for $x$ in the range of $X$ and such that the restriction of $K_n$ 
 to the kernel of $X$ is $-nI_d$, where $d \geq 1$ is the dimension of this kernel.  We have $\tr (K_nX) - \Lambda(K_n) = a d \log (n+1)$ and
 the supremum in (\ref{lambdastar}) is infinite.

If $\det X \not=0$, make the variable change $K=I_p -aX^{-1}L$ and observe that \begin{equation}\label{hineq}
\log \det L  \leq \tr (L -I_p)
\end{equation}
 with equality only at $L=I_p$.
\hfill $\Box$

At last, we have another LDP for rescaled Wishart distributions.  Its proof is left to the reader and uses directly the density (\ref{defwish}).

\begin{prop}
\label{newlem}
Let $p$ and $a$ be fixed. 
If $X$ is $W_p(a)$ distributed then $X/n$ satisfies a LDP in ${\mathcal S}_p^+(\C)$ with good rate function
\begin{equation}
\label{newrate}
{\mathcal I}_s (X) =  \tr X .
\end{equation}
\end{prop}
\bigskip

\noindent{\bf Acknowledgements.}
The authors would like to thank two anonymous referees for their constructive
 comments on an earlier version of this paper.\\
The work of the authors was supported by the Deutsche Forschungsgemeinschaft:
(Sonderforschungsbereich Tr/12; project C2, Fluctuations and universality of invariant random matrix ensembles).
A.R.'s work was partly supported by the ANR project Grandes Matrices
Al\'eatoires ANR-08-BLAN-0311-01.

\bibliographystyle{apalike}
\bibliography{ldpbib3}

\begin{thebibliography}{}

\bibitem[Anderson et~al., 2010]{agz}
Anderson, G., Guionnet, A., and Zeitouni, O. (2010).
\newblock {\em An introduction to random matrices}.
\newblock Cambridge University Press, Cambridge.

\bibitem[Berg, 2008]{moren2008coimbra}
Berg, C. (2008).
\newblock The matrix moment problem.
\newblock In Moren, A. and Branquinho, A., editors, {\em Coimbra Lecture Notes
  on Orthogonal Polynomials}, pages 1--56. Nova Science Pub Inc.

\bibitem[Cartan, 1967]{oldcar}
Cartan, H. (1967).
\newblock {\em Calcul diff\'erentiel}.
\newblock Hermann, Paris.

\bibitem[Chang et~al., 1993]{Changetal}
Chang, F., Kempermann, J., and Studden, W. (1993).
\newblock A normal limit theorem for moment sequences.
\newblock {\em Annals of Probability}, 21(3):1295--1309.

\bibitem[Collins, 2005]{collins2005product}
Collins, B. (2005).
\newblock {Product of random projections, {J}acobi ensembles and universality
  problems arising from free probability}.
\newblock {\em Probability theory and related fields}, 133(3):315--344.

\bibitem[Damanik et~al., 2008]{damanik2008analytic}
Damanik, D., Pushnitski, A., and Simon, B. (2008).
\newblock {The analytic theory of matrix orthogonal polynomials}.
\newblock {\em Surveys in Approximation Theory}, 4:1--85.

\bibitem[Dembo and Zeitouni, 1998]{demboz98}
Dembo, A. and Zeitouni, O. (1998).
\newblock {\em Large Deviations Techniques and Applications}.
\newblock Springer.

\bibitem[Dette and Gamboa, 2007]{DeGa}
Dette, H. and Gamboa, F. (2007).
\newblock Asymptotic properties of the algebraic moment range process.
\newblock {\em Acta Math. Hungar.}, 116(3):247--264.

\bibitem[Dette and Nagel, 2010]{denag09}
Dette, H. and Nagel, J. (2010).
\newblock Matrix measures, random moments and {G}aussian ensembles.
\newblock {\em J. Theor. Prob.}, DOI: 10.1007/s10959-011-0370-7.

\bibitem[Dette and Studden, 1997]{DeSt97}
Dette, H. and Studden, W. (1997).
\newblock {\em The theory of canonical moments with applications in statistics,
  probability, and analysis}.
\newblock Wiley Series in Probability and Statistics,.

\bibitem[Dette and Studden, 2002]{destu02}
Dette, H. and Studden, W.~J. (2002).
\newblock Matrix measures, moment spaces and {F}avard's theorem for the
  interval [0,1] and [0,$\infty$).
\newblock {\em Linear Algebra and its Applications}, 345:169--193.

\bibitem[Dette and Wagener, 2010]{dewag09}
Dette, H. and Wagener, J. (2010).
\newblock Matrix measures on the unit circle, moment spaces, orthogonal
  polynomials and the {G}eronimus relations.
\newblock {\em Linear Algebra and its Applications}, 432:1609--1626.

\bibitem[Dubovoj et~al., 1992]{dubovoj1992matricial}
Dubovoj, V., Fritzsche, B., and Kirstein, B. (1992).
\newblock {\em {Matricial version of the classical Schur problem}}.
\newblock BG Teubner Gmbh.

\bibitem[Fischmann et~al., 2011]{fischmann1}
Fischmann, J., Bruzda, W., Khoruzhenko, B.~A., Sommers, H.-J., and Zyczkowski,
  K. (2011).
\newblock Induced ginibre ensemble of random matrices and quantum operations.
\newblock {\em arXiv.org}, arXiv:1107.5019v1 [math-ph].

\bibitem[Forrester, 2010]{forrester2010}
Forrester, P. (2010).
\newblock {\em Log-Gases and Random Matrices}.
\newblock Princeton University Press.

\bibitem[Forrester and Krishnapur, 2009]{forrester2009derivation}
Forrester, P. and Krishnapur, M. (2009).
\newblock {Derivation of an eigenvalue probability density function relating to
  the {P}oincar\'e disk}.
\newblock {\em Journal of Physics A: Mathematical and Theoretical}, 42:385204.

\bibitem[Gamboa and Lozada-Chang, 2004]{Gamboa}
Gamboa, F. and Lozada-Chang, L.-V. (2004).
\newblock Large deviations for random power moment problem.
\newblock {\em The Annals of Probability}, 32(3B):2819--2837.

\bibitem[Gamboa and Rouault, 2010]{gamrou2010}
Gamboa, F. and Rouault, A. (2010).
\newblock Canonical moments and random spectral measures.
\newblock {\em J. Theor. Probab.}, DOI 10.1007/s10959-009-0239-1.

\bibitem[Ginibre, 1965]{gigi}
Ginibre, J. (1965).
\newblock Statistical ensembles of complex, quaternion, and real matrices.
\newblock {\em J. Mathematical Phys.}, 6:440--449.

\bibitem[Hiai and Petz, 2006]{HiaiPetz}
Hiai, F. and Petz, D. (2006).
\newblock Large deviations for functions of two random projections.
\newblock {\em Acta Sci. Math. (Szeged)}, 72:581--609.

\bibitem[Horn and Johnson, 1985]{hornj85}
Horn, R.~A. and Johnson, C.~R. (1985).
\newblock {\em Matrix Analysis}.
\newblock Cambridge University Press.

\bibitem[Jiang, 2005]{jiang2005maxima}
Jiang, T. (2005).
\newblock {Maxima of entries of {H}aar distributed matrices}.
\newblock {\em Probability Theory and Related Fields}, 131(1):121--144.

\bibitem[Khatri, 1965]{khatri1965}
Khatri, C.~G. (1965).
\newblock Classical statistical analysis based on a certain multivariate
  complex {G}aussian distribution.
\newblock {\em Annals of Mathematical Statistics}, 36:98--114.

\bibitem[Lozada-Chang, 2005]{LozEJP}
Lozada-Chang, L. (2005).
\newblock Large deviations on moment spaces.
\newblock {\em Electronic J. Probab.}, 10:662--690.

\bibitem[Mathai, 1997]{mathai97}
Mathai, A. (1997).
\newblock {\em Jacobians of Matrix Transformations and Functions of Matrix
  Argument}.
\newblock World Scientific Publ.

\bibitem[Mehta, 2004]{mehta04}
Mehta, M. (2004).
\newblock {\em Random matrices, Pure and Applied Mathematics}.
\newblock Elsevier/Acacemic Press, Amsterdam.

\bibitem[Pillai and Jouris, 1971]{piljou1971}
Pillai, K. C.~S. and Jouris, G.~M. (1971).
\newblock Some distribution problems in the multivariate complex {G}aussian
  case.
\newblock {\em Annals of Mathematical Statistics}, 42:517--525.

\bibitem[Robertson and Rosenberg, 1968]{robertson}
Robertson, J. and Rosenberg, M. (1968).
\newblock {The decomposition of matrix-valued measures}.
\newblock {\em Michigan Math. J}, 15:353--368.

\bibitem[Simon, 2005]{simon05}
Simon, B. (2005).
\newblock {\em Orthogonal polynomials on the unit circle. Part 1: Classical
  theory.}
\newblock {Colloquium Publications. American Mathematical Society 54, Part 1.
  Providence, RI: American Mathematical Society (AMS)}.

\bibitem[Skibinsky, 1969]{Ski}
Skibinsky, M. (1969).
\newblock Some striking properties of binomial and beta moments.
\newblock {\em Ann. Math. Statist.}, 40:1753--1764.

\bibitem[van~der Vaart, 1998]{vava}
van~der Vaart, A.~W. (1998).
\newblock {\em Asymptotic statistics}, volume~3 of {\em Cambridge Series in
  Statistical and Probabilistic Mathematics}.
\newblock Cambridge University Press, Cambridge.

\end{thebibliography}

\end{document}